\documentclass[11pt]{article}
\usepackage[a4paper,margin=28mm]{geometry}
\usepackage[T1]{fontenc}
\usepackage{lmodern}
\usepackage{amsmath,amssymb,amsthm,mathtools}
\usepackage{microtype,xurl,hyperref}
\usepackage{authblk}
\hypersetup{colorlinks=true,linkcolor=blue,citecolor=blue,urlcolor=blue}
\theoremstyle{definition}
\newtheorem*{definition}{Definition}
\theoremstyle{plain}
\numberwithin{equation}{section}
\newtheorem{theorem}{Theorem}[section]
\newtheorem{proposition}[theorem]{Proposition}
\newtheorem{lemma}[theorem]{Lemma}
\newcommand{\HH}{\mathbb H}
\newcommand{\BS}{\mathrm{BS}}
\newcommand{\cusp}{\mathrm{cusp}}
\newcommand{\nc}{\mathrm{nc}}
\newcommand{\res}{\mathrm{res}}
\newcommand{\cont}{\mathrm{cont}}
\newcommand{\Ind}{\operatorname{Ind}}
\title{Borel--Serre Callias Classes and Cuspidal Localization}
\author[1]{Hang Wang\thanks{Email: \texttt{wanghang@math.ecnu.edu.cn}}}
\author[1]{Wenqing Wu\thanks{Email: \texttt{wuwenqingsz@qq.com}}}
\affil[1]{Research Center for Operator Algebras,
School of Mathematical Sciences,
East China Normal University,
Shanghai, China}
\begin{document}
\maketitle
\begin{abstract}
We construct equivariant Callias classes on the Borel--Serre space of a
Hilbert modular group and prove a Kasparov product formula.
For torsion-free normal congruence subgroups and twists by powers
$m_j\geq2$ of the holomorphic cotangent line bundles, we identify the
image of the maximal higher index in the representation ring of the
finite quotient with the equivariant index of the cuspidal Dolbeault
operator on total forms.
The comparison uses a Fredholm homotopy that separates an invertible
noncuspidal block.
As an application, we compute the $A_5$-equivariant index on the
principal level-two Hilbert modular surface over $\mathbb Q(\sqrt5)$
for all parallel powers $m\geq2$.
The resulting class in $R(A_5)$ is $15$-periodic in $m$ and has
virtual dimension four.
We express this class as an alternating sum of Hilbert cusp-form
representations and deduce dimension identities for the full Hilbert
modular group.
\end{abstract}

\noindent\textbf{2020 Mathematics Subject Classification.}
Primary 58J20; Secondary 19K35, 19K56, 11F41.

\medskip
\noindent\textbf{Keywords.}
Callias-type operators, Borel--Serre compactification,
equivariant index theory, cuspidal localization, KK-theory.

\section*{Introduction}
\label{sec:introduction}

Let $F$ be a totally real number field of degree $r\geq2$,
and let $\Gamma$ be a Hilbert modular group acting on
$X=\HH^r$. The cuspidal subspace of the $L^2$-space of
automorphic sections is characterised by the vanishing
of constant terms along the unipotent directions at the
cusps. For suitable positive twists, the Dolbeault operator
restricted to this subspace is Fredholm. When a finite group acts on the quotient, this
index takes values in its representation ring and
records the equivariant structure of the cuspidal
spectrum. Our aim is to realise this index through a
geometric operator on the noncompact manifold and to
relate it to an index class in the $K$-theory of a group
$C^*$-algebra.

Our work is motivated by an unpublished manuscript
of G.~G. Kasparov, \emph{Automorphic forms and index
theory}, which is not publicly available.
We extend his approach from the upper half-plane
$\HH$ to Hilbert modular spaces $\HH^r$.
The higher-dimensional cusp geometry and the
noncuspidal automorphic spectrum introduce additional
analytic difficulties. To establish the cuspidal
comparison, we adopt a different analytic route
based on the Eisenstein spectral decomposition
and Shimura's eigenvalue estimates. We show the
absence of residual spectrum for
the twists considered here and obtain a uniform
positive lower bound for the square of the Dolbeault
operator on the orthogonal complement of the cuspidal
subspace. This yields the invertibility required
for the equivariant Fredholm homotopy and allows us
to identify the finite-quotient image of our
Borel--Serre Callias higher index with the
cuspidal index.

As an arithmetic application, we compute the
$A_5$-equivariant index for the Hilbert modular
surface over $\mathbb Q(\sqrt5)$ and express it
as an alternating sum of Hilbert cusp-form
representations. Taking invariants yields explicit
dimension identities for the full Hilbert modular
group. These results connect the higher-index
construction in operator $K$-theory with the
arithmetic theory of Hilbert modular forms.

To formulate these results precisely, fix a torsion-free normal congruence subgroup
$\Gamma_0\triangleleft\Gamma$ of finite index, and put
\[
 M=\Gamma_0\backslash X,\qquad
 Q=\Gamma/\Gamma_0.
\]
The finite group $Q$ acts on $M$ and on the descended
bundles and operators. For integers $m_j\geq2$, let
\[
 L_{\boldsymbol m}
 =\boxtimes_{j=1}^rK_{\HH}^{m_j},
\]
where $K_{\HH}$ denotes the holomorphic cotangent line
bundle on $\HH$. We consider the twisted
Dolbeault--Dirac operator
\[
 D_{\boldsymbol m}
 =\bar\partial_{L_{\boldsymbol m}}
  +\bar\partial_{L_{\boldsymbol m}}^*
\]
on $L_{\boldsymbol m}$-valued total forms, graded
by the parity of total degree. Its descent to $M$
is denoted by $D_{\boldsymbol m}^M$, and its
restriction to the cuspidal subspace by
$D_{\boldsymbol m,\cusp}^M$.

The Borel--Serre space $\overline X^{\BS}$ carries
a proper and cocompact $\Gamma$-action
\cite{BorelSerre}, providing the geometric setting
for maximal analytic assembly. An essential point
is that functions in $C_0(\overline X^{\BS})$ may
have nonzero values at the added boundary, so
compactness of the associated $KK$-cycle requires
estimates extending along the cusp ends. We obtain
these estimates using a smooth $\Gamma$-invariant
function $\beta$ with prescribed growth at the cusps
and the Callias-type perturbation
\[
 \Phi=(\partial\beta\wedge)
       +(\partial\beta\wedge)^*,
 \qquad
 B_{\boldsymbol m}=D_{\boldsymbol m}+\Phi,
\]
where $\Phi$ acts on the holomorphic exterior factor.
The bounded transform of $B_{\boldsymbol m}$ defines
an equivariant class
\[
 x_{\beta,\boldsymbol m}
 \in KK_0^\Gamma
   (C_0(\overline X^{\BS}),\mathbb C).
\]

We also identify this class as a Kasparov product.
Writing $\Theta_\beta$ for the class defined by the
potential and $D_{\boldsymbol m}^{0,*}$ for the
twisted Dolbeault operator on antiholomorphic forms,
we prove
\[
 x_{\beta,\boldsymbol m}
 =\Theta_\beta\widehat\otimes_{C_0(X)}
   [D_{\boldsymbol m}^{0,*}].
\]
This description is compatible with descent to the
finite quotient and identifies the descended class
with that of the Callias operator on $M$.

Maximal analytic assembly gives
\[
 \alpha_{\beta,\boldsymbol m}
 =\Ind_\Gamma(x_{\beta,\boldsymbol m})
 \in K_0(C^*_{\max}(\Gamma)).
\]
The quotient homomorphism
$\Gamma\twoheadrightarrow Q$ induces
$\pi_{\max}:C^*_{\max}(\Gamma)\to C^*(Q)$.
Using the naturality of assembly
\cite{HochsLandsman} and the finite-quotient
identification, we establish the index formula
\[
 \boxed{
 \begin{aligned}
 (\pi_{\max})_*\alpha_{\beta,\boldsymbol m}
 &=\Ind_Q\bigl((B_{\boldsymbol m}^M)^+\bigr)\\
 &=\Ind_Q\bigl(
       (D_{\boldsymbol m,\cusp}^M)^+\bigr)
 \end{aligned}}
 \quad\text{in }R(Q).
\]
Thus the finite-quotient image of the higher index
recovers the full $Q$-equivariant cuspidal index,
including its virtual character values.

The second equality is established by comparing
the descended Callias operator with
\[
 S=D+\Phi P,
\]
where $D$ and $\Phi$ denote the descended operators,
and $P$ is the averaging projection over the unipotent
tori in the cusp regions, extended by zero on the
compact core. The operator $S$ agrees with the
Dolbeault operator on the cuspidal subspace and
has an invertible noncuspidal block.
A relative-positivity argument \cite{vanDenDungen}
provides a $Q$-equivariant Fredholm homotopy
between the bounded transforms of $D+\Phi$ and $S$.

We make the equivariant index explicit for
\[
 F=\mathbb Q(\sqrt5),\qquad
 \Gamma=\mathrm{PSL}_2(\mathcal O_F),\qquad
 \Gamma_0=\Gamma(2),
\]
where reduction modulo $2$ gives
$Q\cong A_5$. For parallel powers
$\boldsymbol m=(m,m)$, $m\geq2$, set
\[
 I_m=\Ind_{A_5}
       \bigl((B_{\boldsymbol m}^M)^+\bigr).
\]
Denoting the irreducible representations of $A_5$
by $\mathbf1,\mathbf3,\mathbf3',\mathbf4,\mathbf5$,
we obtain
\[
 I_m=
 \begin{cases}
  4\mathbf1,
    &15\mid m,\\
  2\mathbf1-\mathbf3-\mathbf3'+2\mathbf4,
    &3\mid m,\ 5\nmid m,\\
  3\mathbf1-\mathbf4+\mathbf5,
    &5\mid m,\ 3\nmid m,\\
  \mathbf1-\mathbf3-\mathbf3'
   +\mathbf4+\mathbf5,
    &3\nmid m,\ 5\nmid m.
 \end{cases}
\]
In particular, $I_{m+15}=I_m$, and every $I_m$
has virtual dimension four. The nonidentity
character values follow from the fixed-point
formula for equivariant index pairings
\cite{HochsWang}. The identity value is obtained
from the ordinary index pairing and equals the
Euler characteristic $\chi(M)=4$.

The cuspidal comparison further expresses $I_m$
in terms of Hilbert cusp forms. If
$S_{k_1,k_2}(\Gamma_0)$ denotes the space of
cusp forms of weight $(k_1,k_2)$, with its
natural $A_5$-action, then
\[
\begin{aligned}
 I_m={}&[S_{2m,2m}(\Gamma_0)]
       -[S_{2m+2,2m}(\Gamma_0)]\\
      &-[S_{2m,2m+2}(\Gamma_0)]
       +[S_{2m+2,2m+2}(\Gamma_0)]
       \quad\text{in }R(A_5).
\end{aligned}
\]
Taking $A_5$-invariants gives, in particular,
\[
\begin{aligned}
 &\dim S_{2m,2m}(\Gamma)
 -\dim S_{2m+2,2m}(\Gamma)\\
 &\quad-\dim S_{2m,2m+2}(\Gamma)
 +\dim S_{2m+2,2m+2}(\Gamma)\\
 &=1+\mathbf1_{3\mid m}
      +2\mathbf1_{5\mid m}.
\end{aligned}
\]
Classical dimension formulas for Hilbert modular
forms provide a number-theoretic setting for
such relations \cite{Hirzebruch73}. The
representation-valued identity above retains
the $A_5$-equivariant information and realises
the alternating sum through a higher index
construction, connecting Hilbert modular forms
with equivariant operator $K$-theory and
noncommutative geometry.

The paper is organised as follows.
Section~\ref{sec:callias} constructs the
Borel--Serre Callias class and identifies its
finite-quotient index.
Section~\ref{sec:homotopy} establishes the
cuspidal comparison through spectral estimates
and a Fredholm homotopy.
Section~\ref{sec:a5} computes the
$A_5$-equivariant index and derives its
consequences for Hilbert cusp forms.

\section*{Acknowledgements}

The authors gratefully acknowledge support from the
National Natural Science Foundation of China
(Grant No. 12271165) and the Science and Technology
Commission of Shanghai Municipality
(Grant Nos. 23JC1401900 and 22DZ2229014).
H. Wang is also partially supported by the Fundamental
and Interdisciplinary Disciplines Breakthrough Plan
of the Ministry of Education of China
(Grant No. JYB2025XDXM112).

During the development of this work, the authors used
OpenAI's ChatGPT to explore mathematical approaches
and proof strategies, including suggestions that
contributed to the development of several key arguments.
ChatGPT was also used to assist with the drafting and
revision of the manuscript.

The authors designed and developed the research
project, critically examined, adapted, and verified
the AI-assisted arguments, and take full
responsibility for the mathematical content
of the paper.

\section{The Borel--Serre Callias class and its finite-quotient index}
\label{sec:callias}

Let $F$ be a totally real field of degree $r$, with ring of
integers $\mathcal O_F$.  Write
\[
 \HH=\{z=x+iy\in\mathbb C:y>0\},\qquad X=\HH^r,
\]
and equip $X$ with its product hyperbolic metric.
Let $\Gamma\subset\mathrm{PSL}_2(\mathbb R)^r$ be a Hilbert
modular group for $F$, commensurable with the image of
$\mathrm{PSL}_2(\mathcal O_F)$ under the real embeddings of $F$
and acting on $X$ by coordinatewise M\"obius transformations.
Choose a torsion-free normal congruence subgroup
$\Gamma_0\triangleleft\Gamma$ of finite index and set
\[
 Q=\Gamma/\Gamma_0,\qquad M=\Gamma_0\backslash X,
\]
so that the finite group $Q$ acts on $M$.
We use a superscript $M$ for bundles and operators descended
from $X$ to $M$.

Our aim is to construct a higher index in
$K_0(C^*_{\max}(\Gamma))$, the $K$-theory of the maximal group
$C^*$-algebra, whose image under $\Gamma\twoheadrightarrow Q$
recovers the $Q$-equivariant cuspidal index on $M$.
To control the noncompact cusp ends, we add a growing potential
to the twisted Dolbeault operator, obtaining a Callias operator.

Write $\overline X^{\BS}$ for the space obtained by adjoining
the Borel--Serre boundary to $X$, equipped with the extended
proper and cocompact $\Gamma$-action.  Its quotient
\[
 \overline M=\Gamma_0\backslash\overline X^{\BS}
\]
is the Borel--Serre compactification of $M$, with the projection
$\overline X^{\BS}\to\overline M$ a $\Gamma_0$-covering
\cite{BorelSerre}.

The growing potential provides the compactness required for
the equivariant Callias class with the action of
$C_0(\overline X^{\BS})$ on sections over $X$ by restriction
and multiplication, even when the functions have nonzero
values on the added boundary.
We express this class as the Kasparov product of the
potential class and the Dolbeault class, and identify its image
under the finite quotient with the $Q$-equivariant index of the
descended Callias operator.  	Section~\ref{sec:homotopy} then identifies this
index with the cuspidal index.

Let $K_{\HH}:=T^{*1,0}\HH$ be the canonical line bundle of
classical weight $2$, with fibre at $z$ the space of
complex-linear covectors on $T_z\HH$ and invariant Hermitian
metric $\|dz\|^2=2y^2$ in the frame $dz$.
For integers $m_j\geq2$, write
$K_{\HH}^{m_j}:=K_{\HH}^{\otimes m_j}$ for its $m_j$-fold tensor
power, generated by $(dz)^{\otimes m_j}$.
Equip
$L_{\boldsymbol m}:=\boxtimes_{j=1}^rK_{\HH}^{m_j}$
with the induced Hermitian metric and Chern connection,
where $\boxtimes$ denotes the tensor product of the pullbacks
from the $r$ factors of $X$.  Put
\begin{equation}\label{eq:model}
 E_{\boldsymbol m}=\Lambda^{*,*}T^*X\otimes L_{\boldsymbol m},
 \qquad D_{\boldsymbol m}=\bar\partial_{L_{\boldsymbol m}}
                          +\bar\partial_{L_{\boldsymbol m}}^*,
\end{equation}
with total grading $(-1)^{p+q}$ on forms of bidegree $(p,q)$.
Here $\bar\partial_{L_{\boldsymbol m}}$ is the Dolbeault
differential on $L_{\boldsymbol m}$-valued forms, and its adjoint
is taken with respect to the specified metrics.
Write $D_{\boldsymbol m}^{0,*}$ for the restriction of
$D_{\boldsymbol m}$ to
$\Lambda^{0,*}T^*X\otimes L_{\boldsymbol m}$.

\subsection{The potential class}
\label{subsec:potential}

We construct a $\Gamma$-equivariant potential whose growth
along the cusps gives the compactness required for the
$C_0(\overline X^{\BS})$-action described above.

A cusp of $\Gamma\backslash X$ is an end corresponding to a
boundary component of $\Gamma\backslash\overline X^{\BS}$.
For each cusp, identify a lift of a neighbourhood with a
neighbourhood of the standard cusp at infinity by an isometry,
and fix the resulting coordinates $z_j=x_j+iy_j$, in which
the product metric is
\[
 ds^2=\sum_{j=1}^r\frac{dx_j^2+dy_j^2}{y_j^2},
 \qquad q:=\prod_{j=1}^r y_j.
\]
The coordinate choices are made separately at the different
cusps; in each case $q$ is invariant under the cusp stabilizer
and descends to the corresponding cusp neighbourhood.
Choose $R_0>0$ large enough that the regions $q>R_0$ are
disjoint cusp neighbourhoods with compact complement, which
we call the core.  Extend their height functions to a smooth
function $q:\Gamma\backslash X\to(0,\infty)$ satisfying
$q\leq R_0$ on the core, and use the same notation for its
pullbacks to $M$ and $X$.

Choose $R_1>R_0$ and a smooth function
$\chi\in C^\infty(\mathbb R)$ vanishing on a neighbourhood of
$(-\infty,R_0]$, with $\chi(t)=t$ for $t\geq R_1$.
For $\varepsilon>0$, define
\[
 \beta:=\varepsilon\,\chi(q)^2.
\]
Thus $\beta$ vanishes on the core and near the cusp entrances,
and equals $\varepsilon q^2$ for $q\geq R_1$.
Use the same notation for its pullbacks to $M$ and $X$;
on the cusps, put
$\dot\beta:=d\beta/d\log q$ and
$\ddot\beta:=d^2\beta/d(\log q)^2$.
Since $|\ddot\beta|-2|\dot\beta|$ vanishes for $q\geq R_1$,
we may take $\varepsilon$ sufficiently small that
\begin{equation}\label{eq:normalization}
 C_\beta:=\sup_{\text{cusp points}}
          \bigl(|\ddot\beta|-2|\dot\beta|\bigr)_+<\tfrac12.
\end{equation}
This bound will be used in the noncuspidal estimate of
Section~\ref{subsec:invertibility}.

\begin{definition}[The potential]
On $S_{\mathrm{hol}}:=\Lambda^{*,0}T^*X$, graded by
holomorphic degree, the potential associated with $\beta$ is
\begin{equation}\label{eq:B}
 \Phi:=(\partial\beta\wedge)+(\partial\beta\wedge)^*,
\end{equation}
where $\partial\beta\wedge$ denotes exterior multiplication
by $\partial\beta$ and the adjoint is taken fibrewise.
We also write $\Phi$ for its extension
$\Phi\widehat\otimes1$ to $E_{\boldsymbol m}$.
\end{definition}

\begin{proposition}\label{prop:potential-class}
The potential $\Phi$ gives rise to a $\Gamma$-equivariant
unbounded Kasparov cycle whose bounded transform defines
\begin{equation}\label{eq:potential-class}
 \Theta_\beta:=[C_0(X,S_{\mathrm{hol}}),\Phi]
       \in KK_0^\Gamma
       \bigl(C_0(\overline X^{\BS}),C_0(X)\bigr).
\end{equation}
\end{proposition}

\begin{proof}
We use the self-adjoint case of the unbounded Kasparov
conditions in \cite[Definition~2.3]{vanDenDungen}.
Regard $\Phi$ as a multiplication operator on
$C_0(X,S_{\mathrm{hol}})$ with its maximal domain.
The fibrewise inverses $(\Phi\pm i)^{-1}$ define bounded
adjointable multiplication operators, so $\Phi$ is regular
and self-adjoint.  By construction, it is odd and
$\Gamma$-equivariant and commutes with the
$C_0(\overline X^{\BS})$-action.

In the cusp metric specified above,
\[
 \partial\log q=\sum_{j=1}^r\frac{dz_j}{2iy_j},
 \qquad |\partial\log q|^2=\frac r2,
\]
and hence
\[
 \Phi^2=|\partial\beta|^2I
       =\frac r2(\dot\beta)^2I
       =2r\varepsilon^2q^4I\qquad(q\geq R_1).
\]
Thus $(1+\Phi^2)^{-1/2}$ tends to zero at the
Borel--Serre boundary, and for every
$a\in C_0(\overline X^{\BS})$,
\[
 a(1+\Phi^2)^{-1/2}
      \in C_0(X,\operatorname{End}S_{\mathrm{hol}}).
\]
It is therefore a compact endomorphism of the Hilbert
$C_0(X)$-module $C_0(X,S_{\mathrm{hol}})$, so the bounded
transform defines the claimed class.
\end{proof}

\subsection{The Callias cycle and its product formula}
\label{subsec:cycle}

To realise the Kasparov product of the potential class
$\Theta_\beta$ and the Dolbeault class, we show that the bounded
transform of $B_{\boldsymbol m}=D_{\boldsymbol m}+\Phi$ defines
the corresponding equivariant class over $\overline X^{\BS}$.

\begin{theorem}\label{thm:cycle}
The closure of $B_{\boldsymbol m}$ on $C_c^\infty(X,E_{\boldsymbol m})$
is odd and self-adjoint. Apply functional calculus using the normalizing function $f(t)=t(1+t^2)^{-1/2}$, and one has
\begin{equation}\label{eq:class}
 x_{\beta,\boldsymbol m}:=[f(B_{\boldsymbol m})]
       \in KK_0^\Gamma(C_0(\overline X^{\BS}),\mathbb C).
\end{equation}
This class satisfies
\begin{equation}\label{eq:product}
 x_{\beta,\boldsymbol m}
   =\Theta_\beta\widehat\otimes_{C_0(X)}
      [D_{\boldsymbol m}^{0,*}].
\end{equation}
\end{theorem}

We first obtain estimates for $\Phi^2$ and
the anticommutator
of $D_{\boldsymbol m}$ and $\Phi$ on the cusps, which will be used to establish the compactness
conditions for the $C_0(\overline X^{\BS})$-action and verify
positivity in the product criterion.

\begin{lemma}\label{lem:clifford}
Writing
$\{D_{\boldsymbol m},\Phi\}
:=D_{\boldsymbol m}\Phi+\Phi D_{\boldsymbol m}$
for the anticommutator, we have on each cusp
\begin{equation}\label{eq:potential-estimate}
 \Phi^2=\frac r2(\dot\beta)^2I,\qquad
 \{D_{\boldsymbol m},\Phi\}\text{ is of order zero},\qquad
 \|\{D_{\boldsymbol m},\Phi\}(x)\|_{\mathrm{op}}
    \leq\frac r2\bigl(|\ddot\beta|+|\dot\beta|\bigr).
\end{equation}
For $q\geq R_1$,
\begin{equation}\label{eq:deep-estimate}
 \Phi^2=2r\varepsilon^2q^4I,\qquad
 \|\{D_{\boldsymbol m},\Phi\}\|_{\mathrm{op}}\leq3r\varepsilon q^2.
\end{equation}
The bounds are independent of $\boldsymbol m$.
\end{lemma}

\begin{proof}
The potential square was computed above. In the unitary
cusp frame, let $e_j,\bar e_j$ denote exterior multiplication
by $dz_j/(\sqrt2 y_j)$ and $d\bar z_j/(\sqrt2 y_j)$,
respectively, with their fibrewise adjoints.
Since $e_je_j^*$ acts by the holomorphic degree
$p_j\in\{0,1\}$, differentiating the unitary frame gives
the $\bar\partial$-coefficient
$2y_j\partial_{\bar z_j}-i(m_j+e_je_j^*)$.
Taking its formal adjoint with respect to the hyperbolic
measure, we obtain $D_{\boldsymbol m}=\sum_jD_j$, where
\[
\begin{aligned}
D_j=\frac1{\sqrt2}\Big[
 &\bar e_j\bigl(2y_j\partial_{\bar z_j}
                -i(m_j+e_je_j^*)\bigr)\\
 &+\bar e_j^*\bigl(-2y_j\partial_{z_j}
                +i(m_j+e_je_j^*-1)\bigr)
\Big].
\end{aligned}
\]
Moreover,
\[
 \Phi=-\frac{i\dot\beta}{\sqrt2}
          \sum_k(e_k-e_k^*).
\]
Since $\beta$ depends only on $\log q$, we have
\[
 [2y_j\partial_{\bar z_j},\dot\beta]
 =[-2y_j\partial_{z_j},\dot\beta]
 =i\ddot\beta.
\]
The exterior algebra relations give
\[
\begin{gathered}
 \{\bar e_j,e_k-e_k^*\}
 =\{\bar e_j^*,e_k-e_k^*\}=0,\\
 [e_je_j^*,e_k-e_k^*]
 =\delta_{jk}(e_j+e_j^*).
\end{gathered}
\]
These identities cancel the first-order terms in
$D_{\boldsymbol m}\Phi+\Phi D_{\boldsymbol m}$, leaving
\[
\begin{aligned}
\{D_{\boldsymbol m},\Phi\}
={}&\frac{\ddot\beta}{2}
 \Bigl(\sum_j(\bar e_j+\bar e_j^*)\Bigr)
 \Bigl(\sum_k(e_k-e_k^*)\Bigr)\\
 &+\frac{\dot\beta}{2}
 \sum_j(e_j+e_j^*)(\bar e_j-\bar e_j^*).
\end{aligned}
\]
The two sums in the first term have norm $\sqrt r$,
and each product in the second has norm $1$.
The claimed bounds follow, using
$\dot\beta=2\varepsilon q^2$ and
$\ddot\beta=4\varepsilon q^2$ for $q\geq R_1$.
\end{proof}

The $q^4$ growth of the potential square dominates the $q^2$
anticommutator term.  This yields the following graph bound and
compactness on the quotient.

\begin{lemma}\label{lem:coercivity}
The operator $B_{\boldsymbol m}$ is essentially self-adjoint on
$C_c^\infty(X,E_{\boldsymbol m})$, and its closure satisfies
\begin{equation}\label{eq:coercive}
 B_{\boldsymbol m}^2\geq D_{\boldsymbol m}^2+cq^4-C
\end{equation}
as a quadratic-form inequality, for some $c,C>0$.
The descended operator $B_{\boldsymbol m}^M$ has compact resolvent.
\end{lemma}

\begin{proof}
Since $X$ is complete and $\Phi$ is a smooth self-adjoint
endomorphism, the perturbed Dolbeault--Dirac operator
$B_{\boldsymbol m}$ is essentially self-adjoint on
$C_c^\infty(X,E_{\boldsymbol m})$
\cite[Theorem~2.1 and \S2.6]{BravermanCecchini}.
The same applies on $M$.

The quartic potential growth in Lemma~\ref{lem:clifford}
dominates the anticommutator and, together with the bounded
coefficients outside the deep cusps, gives
\[
 \|B_{\boldsymbol m}v\|^2+C\|v\|^2
 \geq \|D_{\boldsymbol m}v\|^2+c\|q^2v\|^2,
 \qquad v\in C_c^\infty(X,E_{\boldsymbol m}).
\]
Graph norm approximation extends this estimate to the closure
of $B_{\boldsymbol m}$, proving \eqref{eq:coercive}.
The same bound holds on $M$.

Finally, put $u=(B_{\boldsymbol m}^M-i)^{-1}h$.
Self-adjointness gives
$\|B_{\boldsymbol m}^Mu\|^2+\|u\|^2=\|h\|^2$,
so \eqref{eq:coercive} yields
\[
 \|D_{\boldsymbol m}^Mu\|+\|q^2u\|\leq C\|h\|.
\]
Since $q$ is proper on $M$, choose
$\chi_R\in C_c^\infty(M)$ with $0\leq\chi_R\leq1$
and $\chi_R=1$ on $\{q\leq R\}$.
For fixed $R$, local elliptic regularity gives
$\|\chi_Ru\|_{H^1}\leq C_R\|h\|$.
These sections have support in a fixed compact set, so
Rellich's compact embedding into $L^2$ implies that
$\chi_R(B_{\boldsymbol m}^M-i)^{-1}$ is compact.
Moreover,
\[
 \|(1-\chi_R)(B_{\boldsymbol m}^M-i)^{-1}\|
 \leq CR^{-2}\longrightarrow0.
\]
Thus $(B_{\boldsymbol m}^M-i)^{-1}$ is compact.
\end{proof}

To obtain a cycle over $C_0(\overline X^{\BS})$, we must also check
compactness for coefficients extending to the boundary. The estimate \eqref{eq:coercive} controls the cusp tails in this step.

\begin{lemma}\label{lem:bs-local}
For every $a\in C_0(\overline X^{\BS})$, the operators
$a(1+B_{\boldsymbol m}^2)^{-1/2}$ and $[f(B_{\boldsymbol m}),a]$
are compact on $L^2(X,E_{\boldsymbol m})$.
\end{lemma}

\begin{proof}
For $a\in C_c^\infty(\overline X^{\BS})$, we have
$[B_{\boldsymbol m},a]=c(da)$, where $c$ denotes
Dolbeault Clifford multiplication.
On the support of $a$, the unit coordinates stay
in compact sets, so $y_j\asymp q^{1/r}$.
Smoothness in Borel--Serre cusp coordinates gives
\[
 |da|_g\leq C_a(1+q^{1/r}).
\]
Cut off the tails smoothly between heights $R$ and $2R$.
The compactly supported interior parts of $a$ and $c(da)$
give compact operators after composition with
$(1+B_{\boldsymbol m}^2)^{-1/2}$, by elliptic regularity
and Rellich's lemma. Their tail norms are bounded by
\[
 C_a(R^{-2}+R^{1/r-2})\,
      \|q^2(1+B_{\boldsymbol m}^2)^{-1/2}\|
 \longrightarrow0.
\]
Thus
\begin{equation}\label{eq:compact-local}
 a(1+B_{\boldsymbol m}^2)^{-1/2},\qquad
 c(da)(1+B_{\boldsymbol m}^2)^{-1/2}
 \quad\text{are compact}.
\end{equation}

Applying the L\"owner--Heinz inequality \cite{KatoHeinz}
($0\leq A\leq H$ implies $A^\theta\leq H^\theta$
for $0<\theta\leq1$) in the quadratic-form sense
to \eqref{eq:coercive} with $\theta=1/(2r)$ gives
\[
 \|q^{1/r}(1+B_{\boldsymbol m}^2)^{-1/(4r)}\|\leq C,
\]
so $c(da)(1+B_{\boldsymbol m}^2)^{-1/(4r)}$ is bounded.
For $\mu>0$, put $R_\mu=(1+B_{\boldsymbol m}^2+\mu)^{-1}$.
The bounded transform has the strong integral representation
\begin{equation}\label{eq:f-integral}
 f(B_{\boldsymbol m})
 =\frac1\pi\int_0^\infty
       \mu^{-1/2}B_{\boldsymbol m}R_\mu\,d\mu.
\end{equation}
Spectral calculus also gives
\[
 \|c(da)R_\mu\|
 +(1+\mu)^{-1/2}\|c(da)B_{\boldsymbol m}R_\mu\|
 \leq C_a(1+\mu)^{-1+1/(4r)}.
\]
The preceding estimates also show that multiplication by $a$
preserves $\operatorname{Dom}B_{\boldsymbol m}$. Using \eqref{eq:f-integral} and the resolvent commutator identity
\[
 [R_\mu,a]
 =-B_{\boldsymbol m}R_\mu c(da)R_\mu
  -R_\mu c(da)B_{\boldsymbol m}R_\mu,
\]
we obtain
\begin{align*}
 [f(B_{\boldsymbol m}),a]
   &=\frac1\pi\int_0^\infty\mu^{-1/2}
                  [B_{\boldsymbol m}R_\mu,a]\,d\mu,\\
 [B_{\boldsymbol m}R_\mu,a]
   &=(1+\mu)R_\mu c(da)R_\mu
     -B_{\boldsymbol m}R_\mu c(da)B_{\boldsymbol m}R_\mu.
\end{align*}
Both terms in the second line factor through the compact operator
$c(da)(1+B_{\boldsymbol m}^2)^{-1/2}$
in \eqref{eq:compact-local}.
The integrand is bounded in norm by the integrable function
$C_a\mu^{-1/2}(1+\mu)^{-1+1/(4r)}$,
so the integral converges in norm to a compact operator.
Since $C_c^\infty(\overline X^{\BS})$ is uniformly dense in
$C_0(\overline X^{\BS})$, the assertion follows for arbitrary $a$
by operator-norm continuity.
\end{proof}

The preceding lemmas establish the equivariant Kasparov cycle
in \eqref{eq:class}. It remains to identify its class with the
product in \eqref{eq:product}.

\begin{proof}[Proof of Theorem~\ref{thm:cycle}]
We verify the connection and positivity conditions of
\cite[Theorem~2.2]{vanDenDungen} under the graded unitary
identification
\[
 C_0(X,S_{\mathrm{hol}})\widehat\otimes_{C_0(X)}
 L^2(X,\Lambda^{0,*}T^*X\otimes L_{\boldsymbol m})
       \cong L^2(X,E_{\boldsymbol m}).
\]
The Chern connection on $S_{\mathrm{hol}}$ induces
$D_{\boldsymbol m}$ on total forms; the first-factor potential
becomes $\Phi$.

For the connection condition, first restrict the left algebra to
$C_0(X)$. The operators $\Phi$, $D_{\boldsymbol m}^{0,*}$ and
$B_{\boldsymbol m}$ define unbounded Kasparov modules:
their commutators with $C_c^\infty(X)$ are bounded and their
resolvents are locally compact.
For homogeneous $\xi\in C_c^\infty(X,S_{\mathrm{hol}})$,
put $T_\xi u=\xi\widehat\otimes u$. Then
\[
 B_{\boldsymbol m}T_\xi-(-1)^{|\xi|}T_\xi D_{\boldsymbol m}^{0,*}
       =[D_{\boldsymbol m},T_\xi]_{\mathrm{gr}}+\Phi T_\xi.
\]
The right-hand side and its adjoint have bounded, compactly
supported coefficients involving $\nabla^{S_{\mathrm{hol}}}\xi$
and $\Phi\xi$. Thus the off-diagonal block operator
\[
 \begin{pmatrix}
  0&T_\xi\\
  T_\xi^*&0
 \end{pmatrix}
\]
preserves the domain of
$B_{\boldsymbol m}\oplus D_{\boldsymbol m}^{0,*}$
and has a bounded graded commutator with this operator.
These domain and commutator properties imply the bounded
connection condition for such $\xi$
by \cite[Proposition~2.11]{vanDenDungen}.
Uniform density of $C_c^\infty(X,S_{\mathrm{hol}})$ in
$C_0(X,S_{\mathrm{hol}})$ and norm continuity in $\xi$
then extend this condition to all homogeneous sections,
so it also holds for the left action of $C_0(\overline X^{\BS})$.

For positivity, \eqref{eq:coercive} and the bound
$|\Phi|\leq C(1+q^2)$ for some $C>0$ give
$\operatorname{Dom}B_{\boldsymbol m}
 \subseteq\operatorname{Dom}\Phi$.
Lemma~\ref{lem:clifford} gives
$2\Phi^2+\{D_{\boldsymbol m},\Phi\}\geq-CI$.
Graph norm approximation from $C_c^\infty(X,E_{\boldsymbol m})$
therefore yields
\[
 2\operatorname{Re}\langle B_{\boldsymbol m}u,\Phi u\rangle
 \geq-C\|u\|^2,
 \qquad u\in\operatorname{Dom}B_{\boldsymbol m}.
\]
Proposition~3.10 of \cite{vanDenDungen} gives $t>0$ with
\[
 \{f(B_{\boldsymbol m}),f(t\Phi)\}\geq-I.
\]
Positive rescaling of $\Phi$ gives a norm-continuous equivariant
operator homotopy for $\Theta_\beta$. Multiplying the inequality
by $a$ and $a^*$ gives the localized positivity condition with
$\kappa=1$ for every $a\in C_0(\overline X^{\BS})$. The equivariant bounded
product criterion proves \eqref{eq:product}.
\end{proof}

\subsection{The higher index and its finite quotients}
\label{subsec:quotient}

In this section, we define the higher Callias index and compute
its image under the $K$-theory homomorphism induced by the
quotient map $\Gamma\twoheadrightarrow Q$.
We identify this image with the $Q$-equivariant Fredholm index
of the descended Callias operator.

Let $\Ind_\Gamma$ denote the maximal analytic assembly map.
The higher Callias index of $B_{\boldsymbol m}$ is defined by
\begin{equation}\label{eq:alpha}
 \alpha_{\beta,\boldsymbol m}
       :=\Ind_\Gamma(x_{\beta,\boldsymbol m})
          \in K_0(C^*_{\max}(\Gamma)),
\end{equation}
where $C^*_{\max}(\Gamma)$ denotes the full group $C^*$-algebra
of $\Gamma$.
The quotient map $\Gamma\twoheadrightarrow Q$ induces
$\pi_{\max}:C^*_{\max}(\Gamma)\to C^*(Q)$, 
and the induced homomorphism on $K$-theory maps
$\alpha_{\beta,\boldsymbol m}$ to
$(\pi_{\max})_*\alpha_{\beta,\boldsymbol m}
 \in K_0(C^*(Q))\cong R(Q)$.

Since the $\Gamma$-action on $\overline X^{\BS}$ is proper and cocompact
and $\Gamma_0$ is normal, naturality of the maximal assembly map
\cite{ValetteAssembly}
(see also \cite[Appendix~A, Theorem~A.1]{HochsLandsman})
gives the commutative diagram
\begin{equation}\label{eq:quotient-naturality}
\begin{array}{ccc}
 KK_0^\Gamma(C_0(\overline X^{\BS}),\mathbb C)
 &\xrightarrow{\ \Ind_\Gamma\ }&K_0(C^*_{\max}(\Gamma))
 \\[3pt]
 {\scriptstyle V_{\Gamma_0}}\big\downarrow
 &&\big\downarrow{\scriptstyle(\pi_{\max})_*}
 \\[3pt]
 KK_0^Q(C(\overline M),\mathbb C)
 &\xrightarrow{\ \Ind_Q\ }&K_0(C^*(Q))\cong R(Q).
\end{array}
\end{equation}
Here $V_{\Gamma_0}$ is the quotient homomorphism in equivariant
$K$-homology. With counting measures, the group-algebra map
in Hochs--Landsman's formulation is $\pi_{\max}$.
Thus it remains to identify
$V_{\Gamma_0}(x_{\beta,\boldsymbol m})$, which we do next
using the product formula \eqref{eq:product}.

\begin{lemma}\label{lem:quotient-class}
The quotient class is
\begin{equation}\label{eq:quotient-class}
 V_{\Gamma_0}(x_{\beta,\boldsymbol m})
    =[f(B_{\boldsymbol m}^M)]
    \quad\text{in }KK_0^Q(C(\overline M),\mathbb C).
\end{equation}
\end{lemma}

\begin{proof}
Write $N=\Gamma_0$, $Z=\overline X^{\BS}$ and
$d=[D_{\boldsymbol m}^{0,*}]$.
Since the $\Gamma$-action is proper and $N$ is torsion-free,
$N$ acts freely and properly on both $Z$ and $X$.
The invariant bundle and potential descend to give
\[
 \Theta_\beta^M
   =[C_0(M,S_{\mathrm{hol}}^M),\Phi^M]
   \in KK_0^Q(C(\overline M),C_0(M));
\]
the compactness condition follows from the same potential-growth
estimate as in Proposition~\ref{prop:potential-class}.
For the internal Dolbeault class, the result in
\cite[Theorem~3.12 and \S3.2 fn.~21]{HochsLandsman} applies
on the complete manifold $X$ and gives
\[
 V_N^X(d)=[D_{\boldsymbol m}^{M,0,*}].
\]
Here the superscript records the space used in the quotient map.

We verify compatibility with the potential product for the full
quotient-group action. For $Y=Z,X$, the quotient map
identifies $\Gamma\ltimes Y$ with the pullback of
$Q\ltimes(Y/N)$ along $Y\to Y/N$.
It therefore induces a groupoid equivalence and a Morita class
\[
 p_Y\in KK\bigl(C_0(Y/N)\rtimes Q,
                     C_0(Y)\rtimes_{\max}\Gamma\bigr)
\]
by \cite[Theorem~2.8]{MuhlyRenaultWilliams}.
Let
\[
 J_G:KK^G(A,\mathbb C)
       \xrightarrow{\cong}
       KK(A\rtimes_{\max}G,\mathbb C)
\]
be the dual Green--Julg isomorphism for a discrete group $G$
\cite[Theorem~20.2.7(b)]{Blackadar}.
Write $j_G$ for the maximal descent map.
On properly supported representatives, induction by $p_Y$ is
the averaging construction of \cite[\S3.1]{HochsLandsman}:
it gives the averaged Hilbert space, operator and representations.
Consequently, for an equivariant class represented in this way,
\begin{equation}\label{eq:quotient-morita}
 J_Q(V_N^Y y)=p_Y\widehat\otimes J_\Gamma(y).
\end{equation}

Under these Morita equivalences, the descended potential modules
satisfy
\begin{equation}\label{eq:potential-morita}
 p_Z\widehat\otimes j_\Gamma(\Theta_\beta)
   =j_Q(\Theta_\beta^M)\widehat\otimes p_X.
\end{equation}
Indeed, pullback of $S_{\mathrm{hol}}^M$ to $X$ identifies the
associated bundles over the quotient equivalence, and hence the
completed tensor-product modules with their left and right actions. Since $\Phi$ is invariant and commutes
with multiplication by $C_0(Z)$, the identification also
intertwines the bounded transforms of the potential operators.

The maximal descent map preserves Kasparov products
\cite[Theorem~20.6.2]{Blackadar}. Applying
\eqref{eq:product}, \eqref{eq:quotient-morita} and
\eqref{eq:potential-morita}, we obtain
\[
\begin{aligned}
 J_Q\bigl(V_N^Z(x_{\beta,\boldsymbol m})\bigr)
  &=p_Z\widehat\otimes j_\Gamma(\Theta_\beta)
              \widehat\otimes J_\Gamma(d)\\
  &=j_Q(\Theta_\beta^M)\widehat\otimes p_X
              \widehat\otimes J_\Gamma(d)\\
  &=J_Q\bigl(\Theta_\beta^M\widehat\otimes_{C_0(M)}
                    [D_{\boldsymbol m}^{M,0,*}]\bigr).
\end{aligned}
\]
Injectivity of $J_Q$ identifies the quotient class with this
product in $KK_0^Q(C(\overline M),\mathbb C)$.

Finally, the local compactness, connection and positivity arguments
in the proof of Theorem~\ref{thm:cycle} apply on $M$:
the potential and anticommutator estimates descend unchanged,
and compactly supported smooth sections give the connection
condition. Thus
\[
 \Theta_\beta^M\widehat\otimes_{C_0(M)}
        [D_{\boldsymbol m}^{M,0,*}]
      =[f(B_{\boldsymbol m}^M)]
      \quad\text{in }KK_0^Q(C(\overline M),\mathbb C),
\]
which proves \eqref{eq:quotient-class}.
\end{proof}

\begin{theorem}\label{thm:finite-quotient}
Under the homomorphism on $K$-theory induced by
$\Gamma\twoheadrightarrow Q$, the higher Callias index satisfies
\begin{equation}\label{eq:quotient-index}
 (\pi_{\max})_*\alpha_{\beta,\boldsymbol m}
       =\Ind_Q\bigl((B_{\boldsymbol m}^M)^+\bigr)
       \quad\text{in }R(Q).
\end{equation}
\end{theorem}

\begin{proof}
Apply \eqref{eq:quotient-naturality} to
$x_{\beta,\boldsymbol m}$ and use Lemma~\ref{lem:quotient-class}.
Since $Q$ is finite and $\overline M$ is compact, the assembly image
of this quotient cycle is its $Q$-equivariant Fredholm index.
Lemma~\ref{lem:coercivity} gives compact resolvent for
$B_{\boldsymbol m}^M$, yielding \eqref{eq:quotient-index}.
\end{proof}

Section~\ref{sec:homotopy} relates these finite-quotient indices
to cusp forms by projecting the potential onto the cusp
constant terms.
\section{Cuspidal reduction of the Callias index} \label{sec:homotopy}

In this section, we identify the $Q$-equivariant Callias index
$\Ind_Q((B_{\boldsymbol m}^M)^+)$ obtained in Section~\ref{sec:callias}
with the cuspidal Dolbeault index.  Retain $M=\Gamma_0\backslash\HH^r$, $Q=\Gamma/\Gamma_0$, and $m_j\geq2$, and
write $H=L^2(M,E_{\boldsymbol m}^{M})$, $D=D_{\boldsymbol m}^M$, $\Phi=\Phi^M$, and
$B=B_{\boldsymbol m}^M=D+\Phi$. Let $H_{\cusp}$ be the closure in $H$ of the smooth square-integrable
automorphic sections whose constant terms vanish at every cusp, and put
$H_{\nc}=H_{\cusp}^{\perp}$. This orthogonal splitting reduces $D$
and respects the total grading and the $Q$-action; operator subscripts
denote the corresponding restrictions.

On each chosen cusp $\mathfrak a$, let
$T_{\mathfrak a}=\mathbb R^r/L_{\mathfrak a}$ be its unipotent torus,
with $L_{\mathfrak a}$ the cusp translation lattice.
Choose a cusp frame $(\epsilon_A)$ invariant under the unipotent
translations, and write
$u(x,y)=\sum_Au_A(x,y)\epsilon_A(y)$ in this frame.
With $d\bar x=dx/\operatorname{vol}(T_{\mathfrak a})$ denoting
normalized Haar measure, define
\begin{equation}\label{eq:P-average}
 (Pu)(x,y)
 =\sum_A\left(\int_{T_{\mathfrak a}}u_A(x',y)\,d\bar x'\right)
   \epsilon_A(y).
\end{equation}
The bundle metric is translation invariant, so this averaging is
an orthogonal projection.

Extend this orthogonal projection by zero on the compact core,
with the cusp neighborhoods and unipotent bundle identifications
chosen as lifts of those on $\Gamma\backslash\HH^r$.
Since deck transformations identify the corresponding torus fibres
and preserve their normalized Haar measures, we have
\[
 U_gP=PU_g,\qquad g\in Q.
\]
The radial potential $\Phi$ also commutes with $P$ and
vanishes near the cusp entrances. Thus $\Phi P$ is
symmetric and annihilates cusp forms. With the self-adjoint realizations of
Section~\ref{subsec:invertibility},
\begin{equation}\label{eq:comparison}  S:=D+\Phi P=D_{\cusp}\oplus S_{\nc},\qquad
S_{\nc}:=D_{\nc}+(\Phi P)|_{H_{\nc}}.
\end{equation}

Let $H_{\mathrm{disc}}$ and $H_{\cont}$ denote the discrete and
continuous subspaces in the automorphic spectral decomposition of $H$,
respectively, and put
$H_{\res}=H_{\mathrm{disc}}\ominus H_{\cusp}$.
This decomposition \cite[Prop.~2.10.5]{Labesse} gives
\begin{equation}\label{eq:spectral-splitting}  H_{\nc}=H_{\res}\oplus H_{\cont}.
\end{equation}

Section~\ref{subsec:eisenstein} proves that $H_{\res}=0$ for $m_j\geq2$. 	Section~\ref{subsec:invertibility} then uses
the continuous Eisenstein eigenvalues to obtain a uniform positive lower bound for
$D_{\nc}^2$; together with the potential estimate, this proves that $S_{\nc}$ is
invertible. The $Q$-equivariant Fredholm homotopy in Section~\ref{subsec:homotopy}joins the
bounded transforms of $B$ and $S$, giving $\Ind_Q(B^+)=\Ind_Q(D_{\cusp}^+)$ in
$R(Q)$.

\subsection{The Eisenstein realization of the noncuspidal space}
\label{subsec:eisenstein}

A coefficient sector is specified by a multidegree
$(\boldsymbol p,\boldsymbol q)\in\{0,1\}^r\times\{0,1\}^r$,
where $p_j$ and $q_j$ record the holomorphic and antiholomorphic
exterior degrees in the $j$-th factor.
Labesse's norm formula and Shimura's pole exclusion give
the following identification.
\begin{proposition}\label{prop:eisenstein}
For $m_j\geq2$, every coefficient sector of $E_{\boldsymbol m}^{M}$ satisfies
\begin{equation}\label{eq:eisenstein-completion}
 H_{\nc}=H_{\cont},\qquad H_{\res}=0.
\end{equation}
\end{proposition}

\begin{proof}

We work in a fixed coefficient sector of $H$, with its spectral subspaces
restricted accordingly. Labesse's Proposition~3.4.1 \cite{Labesse} expresses
the squared norm of a pseudo-Eisenstein series through its spectral
coefficient on the purely imaginary parameter axis, assuming that the
contour shift has no residue contribution and that the scalar integrands
decay sufficiently rapidly. We verify these hypotheses by establishing
holomorphy and uniform rapid decay on the swept strip.

Let $G=\operatorname{Res}_{F/\mathbb Q}\mathrm{SL}_2$ and let
$\mathsf P$ be a proper minimal parabolic subgroup. Take the smooth seeds
in Labesse's seed space $D_{\mathsf P}$:
\[
 \phi=\sum_{k=1}^{N}h_k(\log q)v_k,\qquad
 h_k\in C_c^\infty(\mathbb R),\qquad q=\prod_jy_j,
\]
where $K=\mathrm{SO}(2)^r$ and the $v_k$ are $K$-finite Levi character
data compatible with the fixed finite level and coefficient sector.
The associated pseudo-Eisenstein series is the absolutely convergent
sum \cite[\S3.1]{Labesse}
\[
 E_\phi(x):=
 \sum_{\gamma\in\mathsf P(\mathbb Q)\backslash G(\mathbb Q)}
 \phi(\gamma x).
\]
Since $G$ has rational rank one, Labesse's decomposition
\cite[Prop.~3.1.1]{Labesse}, applied to this dense class of smooth seeds,
gives
\begin{equation}\label{eq:two-decompositions}
 H_{\nc}=\overline{\operatorname{span}\{E_\phi\}},
\end{equation}
where $\phi$ ranges over the seeds above.

The height Fourier--Mellin transform $\Psi_\phi(\cdot,\lambda)$ (Labesse's
$\Phi(\cdot,\lambda)$) is entire of Paley--Wiener type and uniformly
rapidly decreasing on bounded vertical strips \cite[\S3.1]{Labesse}.
Its purely imaginary parameter axis is
$\Lambda_{\mathsf P}=i\mathbb R\rho_{\mathsf P}$, where
$\rho_{\mathsf P}$ is the half-sum of the positive roots.

The Eisenstein constant-term formula \cite[\S2.8 and Lem.~2.9.2]{Labesse}
involves the intertwining operators $M(w,\lambda)$, where
$w\in\{1,w_0\}$ and $w_0$ is the nontrivial Weyl element, acting by
$w_0\lambda=-\lambda$.
Writing $\langle\cdot,\cdot\rangle_{\mathsf P}$ for Labesse's inner
product on the inducing data, the scalar-product formula
\cite[Prop.~3.1.2]{Labesse} gives, for
$\lambda_0\in\mathbb R_{>1}\rho_{\mathsf P}$,
\begin{equation}\label{eq:scalar-product-contour}
 \|E_\phi\|_H^2
 =\sum_{w\in\{1,w_0\}}\int_{\lambda_0+\Lambda_{\mathsf P}}
   \left\langle
     M(w,\lambda)\Psi_\phi(\lambda),
     \Psi_\phi(-w\overline\lambda)
   \right\rangle_{\mathsf P}\,d\lambda.
\end{equation}

To shift this contour to $\Lambda_{\mathsf P}$, we express the swept
strip in Shimura's parameter $s$.

In multidegree $(\boldsymbol p,\boldsymbol q)$, with
$p_j,q_j\in\{0,1\}$, set
\begin{equation}\label{eq:weights}
 \sigma_j=2(m_j+p_j-q_j)\geq2,\qquad
 \rho_j=\frac{\sigma_j-i\tau_j}{2},\qquad
 \sum_j\tau_j=0,
\end{equation}
where $\boldsymbol\tau$ specifies a compatible cusp character. With
$\rho$ denoting the tuple of the $\rho_j$, we write
$E_{\rho,\xi}(z,s):=E(z,s;\rho,\xi,\Delta)$ for Shimura's Eisenstein
series in his original normalization, with $\Delta$ lifting $\Gamma_0$
and $\xi$ a $\rho$-regular cusp
\cite[(1.14)--(1.17), (3.10), (4.1), (4.4)]{Shimura}.

With $H_{\mathsf P}$ denoting Labesse's Iwasawa height map and
$q=e^{\langle2\rho_{\mathsf P},H_{\mathsf P}(x)\rangle}$, comparison of
the transformation laws for the split component $A_{\mathsf P}$ at the fixed finite level,
right $K$-type and cusp character
\cite[\S2.8]{Labesse}\cite[(1.6), (4.1)]{Shimura} gives
\begin{equation}\label{eq:parameter}
 \lambda=(2s-1)\rho_{\mathsf P}.
\end{equation}
For $\lambda_0=(2b-1)\rho_{\mathsf P}$ with $b>1$, the swept strip is
$\tfrac12\leq\operatorname{Re}s\leq b$.

For our even integral weights $\sigma_j\geq2$, Shimura's result
\cite[Thm.~4.2, pp.~18--19]{Shimura} shows that these Eisenstein series
are holomorphic on the swept strip. Taking their constant terms at all
cusps and using \cite[Lem.~2.9.2]{Labesse} gives holomorphy of
$M(w,\lambda)$ on the relevant inducing blocks. Together with the
entireness of $\Psi_\phi$, this proves that the scalar integrands in
\eqref{eq:scalar-product-contour} are holomorphic on this strip.

\begin{samepage}
In the present $\mathbb Q$-rank-one setting, Langlands' estimate\footnote{%
In Langlands' rank-one notation the scalar spectral parameter is
$\zeta=\|\rho_{\mathsf P}\|(2s-1)$. His matrix $M(\zeta)$ contains the
nontrivial Weyl intertwining blocks $M(w,\lambda)$ used here; see
\cite[\S6, pp.~83--84, 90]{Langlands}.}
bounds $M(w,\lambda)$ uniformly on the part of the swept strip with
large $\|\operatorname{Im}\lambda\|$, for each inducing block occurring
in $\phi$. Holomorphy gives boundedness on the remaining compact part.
Combining these bounds with the rapid decay of $\Psi_\phi$ shows that, for every
Weyl element $w$ and $N\geq0$, the scalar integrands in
\eqref{eq:scalar-product-contour} satisfy
\begin{equation}\label{eq:contour-decay}
 \left|\left\langle
   M(w,\lambda)\Psi_\phi(\lambda),
   \Psi_\phi(-w\overline\lambda)
 \right\rangle_{\mathsf P}\right|
 \leq C_{\phi,N,b}(1+|t|)^{-N},\quad
 \lambda=(a+it)\rho_{\mathsf P},\quad 0\leq a\leq2b-1.
\end{equation}
\end{samepage}

\begin{samepage}
Let $F_{\phi,\mathsf P}$ denote the coefficient defined in Labesse's
\S3.4 and $\Pi_{\cont}$ the continuous spectral projection. Unfolding
\cite[Lem.~2.4.1]{Labesse} identifies $F_{\phi,\mathsf P}$ with the
continuous Plancherel coefficient of $E_\phi$. Since the holomorphy and
decay estimates above justify shifting the contour in
\eqref{eq:scalar-product-contour} without residues, Labesse's
Proposition~3.4.1 \cite{Labesse} and the continuous-block norm formula
\cite[Thm.~2.10.3]{Labesse} give
\begin{equation}\label{eq:no-residue-norm}
 \|E_\phi\|_H^2
 =\frac12\int_{\Lambda_{\mathsf P}}
   \|F_{\phi,\mathsf P}(\lambda)\|_{\mathsf P}^2\,d\lambda
 =\|\Pi_{\cont}E_\phi\|_H^2.
\end{equation}
\end{samepage}

Thus $E_\phi\in H_{\cont}$ for every seed in this test class. By
\eqref{eq:two-decompositions}, their closed span is $H_{\nc}$, giving
$H_{\nc}=H_{\cont}$ and $H_{\res}=0$.
\end{proof}

For the Laplacian estimate in Section~\ref{subsec:invertibility}, we write Labesse's spectral
expansion \cite[Thm.~2.10.4(ii)--(iii) and footnote~5]{Labesse} in
Shimura's notation on a fixed coefficient sector, using its Petersson
inner product \cite[(2.9)]{Shimura}. For $u\in H_{\nc}$, set
$\widehat u_{\rho,\xi}(\nu):=
\langle u,E_{\rho,\xi}(\cdot,\tfrac12+i\nu)\rangle_H$
in the $L^2$ spectral sense. With the spectral measure $d\mu_{\rho,\xi}$ incorporating the
Weyl and seed-normalization factors and the change of variable
$\lambda=2i\nu\rho_{\mathsf P}$, the expansion in $H$ is
\begin{equation}\label{eq:shimura-expansion}
 u(z)=\sum_{\rho,\xi}\int_{\mathbb R}
 \widehat u_{\rho,\xi}(\nu)
 E_{\rho,\xi}\left(z,\tfrac12+i\nu\right)
 \,d\mu_{\rho,\xi}(\nu).
\end{equation}
The sum ranges over all compatible cusp characters, encoded by $\rho$
in \eqref{eq:weights}, and their $\rho$-regular cusps $\xi$.

\subsection{A lower bound for the noncuspidal block}
\label{subsec:invertibility}

We retain the holomorphic-degree contribution to the Eisenstein
eigenvalues in order to control the terms of $\{D,\Phi\}$ that
change that degree.

In the $(\boldsymbol p,\boldsymbol q)$-sector of
\eqref{eq:weights}, put
\[
a_j=m_j+p_j,\qquad
\sigma_j=2(a_j-q_j).
\]
Let $D_j$ denote the Dolbeault component along the
$j$-th hyperbolic factor. The graded product structure
gives $D=\sum_jD_j$ and $D^2=\sum_jD_j^2$.

We identify the degree-one component with a
line-bundle-valued scalar function by the unitary
bundle isomorphism
\[
g\,d\bar z_j\otimes(dz_j)^{a_j}
\longmapsto
2y_j^2g\,(dz_j)^{a_j-1}.
\]
Under this identification, a direct local computation gives
\[
D_j^2\big|_{q_j=0}
=\frac12 L_j^{2a_j},
\qquad
D_j^2\big|_{q_j=1}
=\frac12 L_j^{2(a_j-1)}+(a_j-1),
\]
where $L_j^\sigma$ is Shimura's weight-$\sigma$
Laplacian \cite[(2.4)]{Shimura}.
Equivalently, on the corresponding scalar coefficients,
\[
D_j^2=\frac12 L_j^{\sigma_j}
       +q_j(a_j-1).
\]

Unitary conjugation by the weight
$\prod_j y_j^{\sigma_j/2}$ then gives
\begin{equation}\label{eq:dirac-laplacian}
D_j^2\big|_{(\boldsymbol p,\boldsymbol q)}
=\frac12\left[
-y_j^2(\partial_{x_j}^2+\partial_{y_j}^2)
+i\sigma_j y_j\partial_{x_j}
+a_j(a_j-1)
\right].
\end{equation}

At $s=\tfrac12+i\nu$, Shimura's eigenvalue formula
\cite[(4.2)]{Shimura} therefore yields
\[
\begin{aligned}
\lambda(D_j^2)
&=\frac12\left[
(a_j-q_j-\tfrac12)^2
+(\nu+\tau_j/2)^2
\right]+q_j(a_j-1)\\
&=\frac12\left[
(a_j-\tfrac12)^2
+(\nu+\tau_j/2)^2
\right].
\end{aligned}
\]
In particular, the resulting spectral threshold is
independent of the antiholomorphic degree $q_j$.

In the $(\boldsymbol p,\boldsymbol q)$-sector of
\eqref{eq:weights}, put
\[
a_j=m_j+p_j,\qquad
\sigma_j=2(a_j-q_j).
\]
Let $D_j$ denote the Dolbeault component along the
$j$-th hyperbolic factor. The graded product structure
gives $D=\sum_jD_j$ and $D^2=\sum_jD_j^2$.

We identify the degree-one component with a
line-bundle-valued scalar function by the unitary
bundle isomorphism
\[
g\,d\bar z_j\otimes(dz_j)^{a_j}
\longmapsto
2y_j^2g\,(dz_j)^{a_j-1}.
\]
Under this identification, a direct local computation gives
\[
D_j^2\big|_{q_j=0}
=\frac12 L_j^{2a_j},
\qquad
D_j^2\big|_{q_j=1}
=\frac12 L_j^{2(a_j-1)}+(a_j-1),
\]
where $L_j^\sigma$ is Shimura's weight-$\sigma$
Laplacian \cite[(2.4)]{Shimura}.
Equivalently, on the corresponding scalar coefficients,
\[
D_j^2=\frac12 L_j^{\sigma_j}
       +q_j(a_j-1).
\]

Unitary conjugation by the weight
$\prod_j y_j^{\sigma_j/2}$ then gives
\begin{equation}\label{eq:dirac-laplacian}
D_j^2\big|_{(\boldsymbol p,\boldsymbol q)}
=\frac12\left[
-y_j^2(\partial_{x_j}^2+\partial_{y_j}^2)
+i\sigma_j y_j\partial_{x_j}
+a_j(a_j-1)
\right].
\end{equation}

At $s=\tfrac12+i\nu$, Shimura's eigenvalue formula
\cite[(4.2)]{Shimura} therefore yields
\[
\begin{aligned}
\lambda(D_j^2)
&=\frac12\left[
(a_j-q_j-\tfrac12)^2
+(\nu+\tau_j/2)^2
\right]+q_j(a_j-1)\\
&=\frac12\left[
(a_j-\tfrac12)^2
+(\nu+\tau_j/2)^2
\right].
\end{aligned}
\]
In particular, the resulting spectral threshold is
independent of the antiholomorphic degree $q_j$.

For $u\in\operatorname{Dom}D_{\nc}$, write $u_{\boldsymbol p,\boldsymbol q}$
for its multidegree components. The contractions of
Lemma~\ref{lem:clifford} satisfy
$\|e_j^*u\|^2=\sum_{\boldsymbol p,\boldsymbol q}
p_j\|u_{\boldsymbol p,\boldsymbol q}\|^2$.
Since $D_{\nc}^2$ preserves these multidegrees, their orthogonality
and the unitary expansion \eqref{eq:shimura-expansion} give,
after dropping the nonnegative frequency terms,
\begin{equation}\label{eq:gap}
 \begin{aligned}
 \|D_{\nc}u\|^2
 &\geq\frac12\sum_{\boldsymbol p,\boldsymbol q}\sum_j
       (m_j+p_j-\tfrac12)^2
       \|u_{\boldsymbol p,\boldsymbol q}\|^2\\
 &=\frac12\sum_j(m_j-\tfrac12)^2\|u\|^2
       +\sum_jm_j\|e_j^*u\|^2\\
 &\geq\frac{9r}{8}\|u\|^2+2\sum_j\|e_j^*u\|^2,
 \qquad m_j\geq2.
 \end{aligned}
\end{equation}
Using $p_j^2=p_j$, we obtain
$\frac12(m_j+p_j-\tfrac12)^2
=\frac12(m_j-\tfrac12)^2+m_jp_j$,
which gives the equality. The pointwise norm $|e_j^*u|$ is the
norm of the component with $p_j=1$ and is independent of the
local orthonormal frame. This degree contribution in \eqref{eq:gap}
will control the $\dot\beta$-part of the anticommutator.
To apply this estimate, we realize $S$ initially on the domain
$\mathcal C=C_c^\infty(M,E_{\boldsymbol m}^{M})$.

\begin{proposition}\label{prop:invertibility}
For $m_j\geq2$ and $C_\beta<\tfrac12$, the closure of $S$ on
$\mathcal C$ is self-adjoint and has the decomposition
\eqref{eq:comparison}. Its noncuspidal block $S_{\nc}$ has a
bounded inverse.
\end{proposition}

\begin{proof}
The cusp coefficients of $D$ are independent of $x$, and $\Phi$
vanishes near the cusp entrances. Thus $S$ preserves $\mathcal C$,
and
\[
 \{D,\Phi P\}=\{D,\Phi\}P,
\]
with the order-zero coefficient commuting with $P$.
Since $\Phi P$ is bounded on bounded-height regions and commutes
with radial exhaustion cutoff functions, the cutoff argument
in \cite[\S4]{BravermanCecchini} applies and proves essential
self-adjointness of $S$ on $\mathcal C$.
Cutoff approximation, using $\Phi P=0$ on cusp forms, identifies
$D_{\cusp}$ as a self-adjoint restriction of $S$, giving
\eqref{eq:comparison}. Averaging preserves the grading and the
$Q$-action, so $S$ is odd and $Q$-equivariant.

The projection of $\mathcal C$ onto $H_{\nc}$ is a core for
$S_{\nc}$. For $u$ in this core, put $v=Pu$.
Since $P$ is orthogonal and preserves the multidegrees,
$\|e_j^*v\|\leq\|e_j^*u\|$. Combining \eqref{eq:gap} with
$\|\Phi v\|^2=\tfrac r2\|\dot\beta v\|^2$ gives
\[
 \begin{aligned}
 \|S_{\nc}u\|^2
 &=\|D_{\nc}u\|^2+\|\Phi v\|^2
       +\langle v,\{D,\Phi\}v\rangle\\
 &\geq\frac{9r}{8}\|u\|^2+2\sum_j\|e_j^*v\|^2
       +\frac r2\|\dot\beta v\|^2
       +\langle v,\{D,\Phi\}v\rangle.
 \end{aligned}
\]
In the explicit formula of Lemma~\ref{lem:clifford}, the
$\dot\beta$-part of $\{D,\Phi\}$ is
$\tfrac12\dot\beta\sum_j(e_j+e_j^*)(\bar e_j-\bar e_j^*)$.
By the exterior algebra relations, its quadratic form is the
sum of real inner products below. Estimating the
$\ddot\beta$-part by its norm and then applying Cauchy--Schwarz
and Young's inequality gives
\[
 \begin{aligned}
 \langle v,\{D,\Phi\}v\rangle
 &\geq-\frac r2\langle v,|\ddot\beta|v\rangle
       +\sum_j\operatorname{Re}
          \langle e_j^*v,\dot\beta(\bar e_j-\bar e_j^*)v\rangle\\
 &\geq-\frac r2\langle v,|\ddot\beta|v\rangle
       -\sum_j\|e_j^*v\|\,\|\dot\beta v\|\\
 &\geq-\frac r2\langle v,|\ddot\beta|v\rangle
       -2\sum_j\|e_j^*v\|^2-\frac r8\|\dot\beta v\|^2.
 \end{aligned}
\]
Substituting this lower bound cancels the degree terms; completing
the remaining square using \eqref{eq:normalization} gives
\begin{equation}\label{eq:noncusp-invertibility}
 \begin{aligned}
 \|S_{\nc}u\|^2
 &\geq\frac{9r}{8}\|u\|^2
       +2\sum_j\|e_j^*v\|^2+\frac r2\|\dot\beta v\|^2\\
 &\qquad-\frac r2\langle v,|\ddot\beta|v\rangle
       -2\sum_j\|e_j^*v\|^2-\frac r8\|\dot\beta v\|^2\\
 &=\frac{9r}{8}\|u\|^2
       +\frac{3r}{8}\|\dot\beta v\|^2
       -\frac r2\langle v,|\ddot\beta|v\rangle\\
 &\geq\frac{9r}{8}\|u\|^2
       +\frac{3r}{8}\|(|\dot\beta|-\tfrac43)v\|^2
       -r\left(\frac23+\frac{C_\beta}{2}\right)\|v\|^2\\
 &\geq r\left(\frac{11}{24}-\frac{C_\beta}{2}\right)\|u\|^2
 \geq\frac{5r}{24}\|u\|^2.
 \end{aligned}
\end{equation}
The last line uses $\|v\|\leq\|u\|$ and $C_\beta<\tfrac12$.
Core approximation extends this bound to the domain of $S_{\nc}$;
self-adjointness then gives a bounded inverse.
\end{proof}
\subsection{Comparison of the bounded transforms}
\label{subsec:homotopy}

The invertible noncuspidal block of $S$ has index zero.
We transfer this reduction to $B$ by a Fredholm homotopy of
their bounded transforms. The domain and compactness properties
below allow us to apply relative positivity and work in the
Calkin algebra.

Write $f(t)=t(1+t^2)^{-1/2}$ and let
$\pi:\mathcal B(H)\to\mathcal B(H)/\mathcal K(H)$ be the Calkin quotient.

\begin{lemma}\label{lem:domains}
The operators $B$ and $S$ have the common invariant core $\mathcal C$.
Moreover, $S$ is defined on the domain of $B$ and has compact resolvent.
\end{lemma}

\begin{proof}
Proposition~\ref{prop:invertibility} and Lemma~\ref{lem:coercivity}
give the common core. The graph bound in the latter lemma, together
with $\|\Phi Pu\|\leq\|\Phi u\|$, gives
\[
 \|Su\|\leq\|Du\|+\|\Phi u\|
       \leq C\bigl(\|Bu\|+\|u\|\bigr),\qquad u\in\mathcal C.
\]
Core approximation gives the domain inclusion.

For compactness, Lemma~\ref{lem:clifford} and
\eqref{eq:normalization} give, by completing the square,
\begin{equation}\label{eq:S-coercivity}
 \|Su\|^2
 \geq\|Du\|^2+\tfrac12\|\Phi Pu\|^2
       -r\left(\tfrac94+\tfrac{C_\beta}{2}\right)\|u\|^2,
 \qquad u\in\mathcal C.
\end{equation}
Applied to differences of core approximants, this bound extends
to the domain of $S$ by closedness of $D$ and $\Phi P$.

For a nonzero Fourier frequency $\eta$, write $\eta_j$ for its
real embeddings. On each cusp the frequencies lie in a fixed
fractional ideal, so $\prod_j|\eta_j|$ is bounded below by a
positive constant. Uniformly over the finitely many cusps,
\[
 \sum_jy_j^2|\eta_j|^2
 \geq r\bigl(q^2\prod_j|\eta_j|^2\bigr)^{1/r}
 \geq c q^{2/r}.
\]
In \eqref{eq:dirac-laplacian}, the $y_j$ derivative forms are
nonnegative; absorbing the linear $x_j$ terms gives
$D^2\geq cq^{2/r}-C$ on these modes in the deep cusp.
On the zero mode, \eqref{eq:deep-estimate}
gives $\Phi^2=2r\varepsilon^2q^4I$.
Together with \eqref{eq:S-coercivity}, these bounds give, after applying fixed
radial cusp cutoff functions with bounded commutators with $S$,
\[
 \int_{\{q\geq R\}}|u|^2
 \leq CR^{-2/r}\bigl(\|Su\|^2+\|u\|^2\bigr),\qquad R\gg1,
\]
first on $\mathcal C$ and then by closure.
Local elliptic regularity and \eqref{eq:S-coercivity} give
Rellich compactness on bounded-height regions. Together with
the uniform tail estimate, this makes the graph-norm unit
ball of $S$ relatively compact in $H$, proving compactness
of the resolvent.
\end{proof}

An estimate for $\{B,S\}$ gives a lower bound for the square
of the connecting path in the Calkin algebra.

\begin{theorem}\label{thm:homotopy}
The bounded transforms of $B$ and $S$ are joined by the
$Q$-equivariant Fredholm homotopy
\begin{equation}\label{eq:bounded-path}
 F_t=(1-t)f(B)+tf(S),\qquad 0\leq t\leq1,
\end{equation}
and
\begin{equation}\label{eq:index-equality}
 \Ind_Q(B^+)=\Ind_Q(D_{\cusp}^+)
 \qquad\text{in }R(Q).
\end{equation}
\end{theorem}

\begin{proof}
On $\mathcal C$,
\[
 \{B,S\}=2D^2+2\Phi^2P+R,\qquad R=\{D,\Phi\}(1+P).
\]
Since $P$ commutes with $\Phi$ and $\{D,\Phi\}$, the first two
terms are positive and $R$ is symmetric.
Lemma~\ref{lem:clifford}, \eqref{eq:normalization}, and
Lemma~\ref{lem:coercivity} give
\[
 \|Ru\|\leq C(\|\Phi u\|+\|u\|)
          \leq C'(\|Bu\|+\|u\|).
\]
Core approximation extends $R$ symmetrically to the domain of $B$.
The common invariant core and the domain inclusion in
Lemma~\ref{lem:domains} verify the remaining hypotheses of
\cite[Lemma~3.5, pp.~12--13]{vanDenDungen}. Applying that lemma to $(B,S)$ gives
\[
 2\operatorname{Re}\langle Bu,Su\rangle
 \geq-c\langle u,(1+B^2)^{1/2}u\rangle,
 \qquad u\text{ in the domain of }B.
\]
By \cite[Prop.~3.10, pp.~15--16]{vanDenDungen}, there is $a>0$
such that $\{f(B),f(aS)\}\geq-I$.
Since $f(a\,\cdot)-f$ belongs to $C_0(\mathbb R)$, compact resolvent of $S$
gives $f(aS)-f(S)\in\mathcal K(H)$, hence
\[
 \pi(\{f(B),f(S)\})\geq-1.
\]

Compact resolvent of $B$ and $S$
(Lemma~\ref{lem:coercivity} and
Lemma~\ref{lem:domains}) gives $\pi(f(B))^2=\pi(f(S))^2=1$.
Together with the anticommutator bound, squaring
\eqref{eq:bounded-path} yields
\[
 \pi(F_t)^2\geq\bigl(1-3t(1-t)\bigr)1\geq\tfrac14\,1.
\]
Thus the path is Fredholm and preserves the index in $R(Q)$.
By \eqref{eq:comparison},
$f(S)=f(D_{\cusp})\oplus f(S_{\nc})$; the second summand is invertible
by Proposition~\ref{prop:invertibility} and has graded index zero.
Additivity proves \eqref{eq:index-equality}.
\end{proof}

Together with Theorem~\ref{thm:finite-quotient} and
\eqref{eq:quotient-index}, this identifies the image of
$\alpha_{\beta,\boldsymbol m}$ under $(\pi_{\max})_*$
with the $Q$-equivariant cuspidal Dolbeault index:
\begin{equation}\label{eq:higher-cuspidal-index}
 (\pi_{\max})_*\alpha_{\beta,\boldsymbol m}
   =\Ind_Q(D_{\cusp}^+).
\end{equation}
\section{An explicit equivariant index on a Hilbert modular surface}
\label{sec:a5}

We compute the finite-quotient index of Section~\ref{sec:callias}
for the principal level-two Hilbert modular surface over
$\mathbb Q(\sqrt5)$. Its nonidentity character values follow from
the compactly supported index pairing, while an Euler-class
calculation determines the value at the identity. The resulting
virtual representation gives an identity among four spaces of
Hilbert cusp forms.

Throughout this section, put
\[
 F=\mathbb Q(\sqrt5),\qquad \mathcal O_F=\mathbb Z[u],\qquad
 u=\frac{1+\sqrt5}{2},\qquad u'=1-u=-u^{-1},
\]
and let
\[
 \Gamma=\operatorname{PSL}_2(\mathcal O_F),\qquad
 \Gamma_0=\Gamma(2),\qquad M=\Gamma_0\backslash\HH^2.
\]
Here $\Gamma(2)$ is the kernel of reduction modulo $2$ in
$\operatorname{PSL}_2(\mathcal O_F)$. Elementary matrices give a
surjective reduction map, and hence
\begin{equation}\label{eq:a5-quotient}
 Q=\Gamma/\Gamma_0\cong\operatorname{SL}_2(\mathbb F_4)\cong A_5.
\end{equation}

Let $K_j$ be the holomorphic cotangent line along the $j$-th
factor, so $K_M=K_1\otimes K_2$. For parallel weights \(\boldsymbol m=(m,m)\), \(m\geq2\), write
\[
E_m=E_{\boldsymbol m}^M,\qquad
D_m=D_{\boldsymbol m}^M,\qquad
B_m=B_{\boldsymbol m}^M.
\]
Define
\begin{equation}\label{eq:a5-index}
 I_m=\Ind_Q(B_m^+)\in R(Q),\qquad
 \chi_m(g)=\operatorname{Tr}(g\mid I_m),
\end{equation}
where the trace on a virtual representation means the difference
of the two traces.

\subsection{The equivariant fixed-point formula}
\label{subsec:a5-formula}

The product formula of Section~\ref{subsec:cycle} descends to a
compactly supported pairing on $M$. This identifies the Fredholm
index to which the noncompact fixed-point theorem will be applied.

\begin{proposition} \label{prop:a5-pairing}
The potential defines a compactly supported class
\[
\Theta_{\beta,c}
=[C_0(M,\Lambda^{*,0}T^*M),f(\Phi)]
\in KK_Q^0(\mathbb C,C_0(M)).
\]
Moreover,
\begin{equation}\label{eq:a5-pairing}
I_m=\Theta_{\beta,c}
\widehat\otimes_{C_0(M)}[D_m^{0,*}]
=\operatorname{Ind}_Q(D_{m,\mathrm{cusp}}^+).
\end{equation}
\end{proposition}
\begin{proof} Since \(\Phi^2\) tends to infinity along every cusp, its bounded transform represents a compactly supported class. This is the descended potential class of 	Section~\ref{subsec:quotient} with its left action restricted to the constants. The product formula of Theorem~\ref{thm:cycle} and its descent in Lemma~\ref{lem:quotient-class} identify the pairing with \(\operatorname{Ind}_Q(B_m^+)\). Theorem~\ref{thm:homotopy} then gives the cuspidal index equality.
\end{proof}

To compute the nonidentity character values, we use
Hochs--Wang's fixed-point formula for index pairings
\cite[Theorem~2.18]{HochsWang}. It expresses the
character of the pairing between a compactly supported
equivariant $K$-theory class and the $K$-homology class
of an odd, essentially self-adjoint elliptic operator
in terms of data on the fixed-point set. The theorem
requires a compact Abelian group generated densely
by the chosen element, an isometric action preserving
the operator, and a compact fixed-point set.

For $g\in Q\setminus\{e\}$, take the finite cyclic group
$C=\langle g\rangle$. Its action preserves the complete
hyperbolic metric, the holomorphic structure and the
Chern connection on $M$, so $D_m^{0,*}$ satisfies the
operator hypotheses. Proposition~3.1 supplies the
compactly supported class $\Theta_{\beta,c}$, while
Lemma~\ref{lem:a5-fixed-count} shows that $M^g$ is finite. Restricting
the pairing of Proposition~\ref{prop:a5-pairing} to $C$ and evaluating
at $g$, the fixed-point formula gives
\begin{equation}\label{eq:fixed-point-pairing2}
\chi_m(g)=
\int_{TM^g}
\frac{
\operatorname{ch}(\Theta_{\beta,c}|_{M^g}(g))
\operatorname{ch}([\sigma(D_m^{0,*})|_{TM^g}](g))
\operatorname{Todd}(TM^g\otimes\mathbb C)
}{
\operatorname{ch}(\Lambda_{-1}(N_g\otimes\mathbb C)(g))
},
\end{equation}
where $N_g$ is the real normal bundle of $M^g$ in $M$
and $\Lambda_{-1}V=\sum_k(-1)^k\Lambda^kV$.

\begin{proposition}\label{prop:a5-local}
For every $g\ne e$,
\begin{equation}\label{eq:a5-local}
 \chi_m(g)=\sum_{x\in M^g}
    \det\bigl(dg_x|_{T_x^{1,0}M}\bigr)^{-m}.
\end{equation}
\end{proposition}

\begin{proof}
At $x\in M^g$, let $\lambda_1,\lambda_2$ be the eigenvalues of
$A_x=dg_x|_{T_x^{1,0}M}$. They have modulus one and differ from
$1$. The geometric lift to cotangent fibres is
$(dg_x^{-1})^*$; its holomorphic eigenvalues are $\lambda_j^{-1}$,
and its antiholomorphic eigenvalues are
$\overline{\lambda_j}^{-1}=\lambda_j$.
Consequently, $g$ acts on $(K_M^m)_x$ by
$(\lambda_1\lambda_2)^{-m}$.

Invariance of $\beta$ and the absence of invariant tangent
vectors imply $\partial\beta(x)=0$. The restricted potential
class is therefore $\Lambda_{-1}T_x^{*1,0}M$, with character
$\prod_j(1-\lambda_j^{-1})$. Restriction of the Dolbeault
symbol to the zero tangent fibre gives
$\Lambda_{-1}T_x^{*0,1}M\otimes(K_M^m)_x$, with character
$\prod_j(1-\lambda_j)(\lambda_1\lambda_2)^{-m}$.
Finally, $N_{g,x}=T_xM$, so its complexified Euler factor is
$\prod_j(1-\lambda_j)(1-\lambda_j^{-1})$.
The local term of \eqref{eq:fixed-point-pairing2} is thus
\[
 \frac{\prod_{j=1}^2(1-\lambda_j^{-1})
       \prod_{j=1}^2(1-\lambda_j)
       (\lambda_1\lambda_2)^{-m}}
      {\prod_{j=1}^2(1-\lambda_j)(1-\lambda_j^{-1})}
       =\det(A_x)^{-m}.
\]
All denominator factors are nonzero, and the Todd factor at a
point is $1$. Summation proves the assertion.
\end{proof}

\subsection{Fixed points and local contributions}
\label{subsec:a5-calculation}

A point of $M$ fixed by $g\in Q\setminus\{e\}$
lifts to a point of $\HH^2$ fixed by an elliptic
element of $\Gamma$ mapping to $g$.
Thus all nonidentity fixed points lie above
elliptic $\Gamma$-orbits.
Hirzebruch--Van de Ven
\cite[Proposition~II.1]{HirzebruchVanDeVen}
give exactly two elliptic orbits of each order
$2$, $3$, and $5$, and none of order greater
than $5$. Order $4$ is impossible, since an
element of that projective order would have
trace $\pm\sqrt{2}\notin F$.
Their local quotient description also shows
that the point stabilizers are cyclic.
Hence every nontrivial point stabilizer in
$\Gamma$ has order $2$, $3$, or $5$.

To construct representatives of the two
orbits of each order, we use the trace to
prescribe the projective order.
Varying the off-diagonal entries by reciprocal
units preserves the trace and determinant
while allowing different $\Gamma$-orbits.
Since the unit $u$ has norm $-1$, its two
real embeddings have opposite signs.
This motivates the choices $c=1,u$ in
\begin{equation}\label{eq:a5-generators}
t_2=0,\qquad t_3=1,\qquad t_5=u,\qquad
G_{n,c}=
\begin{pmatrix}
0&-c^{-1}\\
c&t_n
\end{pmatrix},
\quad c\in\{1,u\}.
\end{equation}
These matrices are elliptic at both real
embeddings and have projective order $n$.
Write $p_{n,c}\in\HH^2$ for their unique
fixed points and $g_{n,c}\in Q$ for their
images under $\Gamma\to Q$.

The images $g_{n,1}$ and $g_{n,u}$ belong
to the same conjugacy class in $Q$.
Indeed, writing $\bar u\in\mathbb F_4$ for
the reduction of $u$, we have $\bar u^3=1$ and
\[
g_{n,u}
=\operatorname{diag}(\bar u,\bar u^2)\,
g_{n,1}\,
\operatorname{diag}(\bar u,\bar u^2)^{-1}.
\]

\begin{lemma}\label{lem:a5-orbits}
For each $n\in\{2,3,5\}$, the points
$p_{n,1}$ and $p_{n,u}$ represent the two
elliptic $\Gamma$-orbits of order $n$.
Moreover,
\[
\Gamma_{p_{n,c}}=\langle G_{n,c}\rangle,
\qquad c\in\{1,u\}.
\]
\end{lemma}

\begin{proof}
Write
\[
t_{n,j}=\iota_j(t_n),\qquad
\Delta_{n,j}=\sqrt{4-t_{n,j}^2}
\quad (j=1,2).
\]
Solving the fixed-point equation for
$G_{n,c}$ in each real embedding, with
the root chosen in the upper half-plane,
gives
\[
p_{n,1}=
\left(
\frac{-t_{n,1}+i\Delta_{n,1}}{2},
\frac{-t_{n,2}+i\Delta_{n,2}}{2}
\right)
\]
and
\[
p_{n,u}=
\left(
\frac{-t_{n,1}+i\Delta_{n,1}}{2u},
\frac{-t_{n,2}-i\Delta_{n,2}}{2u'}
\right).
\]
The signs of the imaginary terms in
$p_{n,u}$ reflect $u>0$ and $u'<0$.

Suppose that
\[
A=\begin{pmatrix}*&*\\a&b\end{pmatrix}
\in\mathrm{SL}_2(\mathcal O_F)
\]
maps $p_{n,1}$ to $p_{n,u}$.
At the first real embedding, the
imaginary-part formula for the Mobius
action gives
\[
u=|a(p_{n,1})_1+b|^2
=a^2-t_nab+b^2.
\]
By injectivity of the first embedding,
this equality holds in $F$. Applying
the second embedding yields
\[
u'=(a')^2-t_n'a'b'+(b')^2.
\]
On the other hand, the imaginary-part
formula at the second embedding, where
$u'<0$, gives
\[
-u'=|a'(p_{n,1})_2+b'|^2
=(a')^2-t_n'a'b'+(b')^2.
\]
Thus $u'=-u'$, a contradiction.
The two points therefore represent
distinct $\Gamma$-orbits, which exhaust
the two orbits of order $n$ by the
classification above.

Since $G_{n,c}$ fixes $p_{n,c}$ and
has projective order $n$, the same
classification gives
$\Gamma_{p_{n,c}}=\langle G_{n,c}\rangle$.
\end{proof}

An elliptic $\Gamma$-orbit gives a finite
$Q$-orbit on $M$, whose points need not
all be fixed by a given $g\in Q$.
We determine the fixed points through
the corresponding coset space.

\begin{lemma}\label{lem:a5-fixed-count}
Let $g\in Q$ have order $n>1$.
Above each elliptic $\Gamma$-orbit of
order $n$, exactly two points of $M$
are fixed by $g$. Consequently, $M^g$
consists of four points.
\end{lemma}

\begin{proof}
We first show that $\Gamma_0$ is torsion-free.
By Lemma~\ref{lem:a5-orbits}, every nontrivial
finite-order element of $\Gamma$ is conjugate
to a nonidentity power of some $G_{n,c}$.
The reduction of $G_{n,c}$ modulo $2$ is
nontrivial, since its upper-left entry
is zero.
As $n=2,3,5$ is prime, its reduction
still has order $n$.
Thus every nonidentity power of
$G_{n,c}$, as well as every conjugate
of such a power, remains nontrivial
in $Q$. Hence $\Gamma_0$ is torsion-free.

For each $c\in\{1,u\}$, Lemma~\ref{lem:a5-orbits}
gives the stabilizer
\[
Q_{[p_{n,c}]}
=\operatorname{im}(\Gamma_{p_{n,c}}\to Q)
=\langle g_{n,c}\rangle,
\]
which is cyclic of order $n$.
The points of $M$ above $\Gamma p_{n,c}$
form the $Q$-set
\[
\Gamma_0\backslash\Gamma/\Gamma_{p_{n,c}}
\cong Q/Q_{[p_{n,c}]}.
\]
A coset $aQ_{[p_{n,c}]}$ is fixed by $g$
precisely when
\begin{equation}\label{eq:a5-coset-count}
gaQ_{[p_{n,c}]}=aQ_{[p_{n,c}]}
\iff
a^{-1}\langle g\rangle a=Q_{[p_{n,c}]}.
\end{equation}
All cyclic subgroups of order $n$ in $A_5$
are conjugate, so we may choose $a_0\in Q$
with $a_0^{-1}\langle g\rangle a_0=Q_{[p_{n,c}]}$.
Writing $a=a_0x$, we obtain
\[
\begin{aligned}
a^{-1}\langle g\rangle a=Q_{[p_{n,c}]}
&\iff x^{-1}Q_{[p_{n,c}]}x=Q_{[p_{n,c}]}\\
&\iff x\in N_Q(Q_{[p_{n,c}]}).
\end{aligned}
\]
Thus the distinct fixed cosets are
indexed by
$N_Q(Q_{[p_{n,c}]})/Q_{[p_{n,c}]}$,
and their number is
\[
\frac{|N_Q(Q_{[p_{n,c}]})|}
{|Q_{[p_{n,c}]}|}
=\frac{|N_Q(Q_{[p_{n,c}]})|}{n}=2,
\]
since the normalizer orders in $A_5$
are $4$, $6$, and $10$ for
$n=2,3,5$, respectively.

Finally, any $g$-fixed point $[z]\in M$
has a lift $z\in\HH^2$ fixed by an
element of $\Gamma$ mapping to $g$.
Since $\Gamma_0$ is torsion-free,
the finite stabilizer $\Gamma_z$
injects into $Q$ and contains an
element of order $n$.
The classification of elliptic
stabilizers gives $|\Gamma_z|=n$.
Lemma~\ref{lem:a5-orbits} accounts
for both elliptic orbits of this
order, so there are no further
fixed points.
\end{proof}

For $c\in\{1,u\}$, put
\[
d_{n,c}=
\det\bigl(
dG_{n,c}|_{T^{1,0}_{p_{n,c}}\HH^2}
\bigr).
\]
The Mobius derivative at each of the
two representatives gives
\begin{equation}\label{eq:a5-derivative}
d_{n,c}
=\prod_{j=1}^2
\left(
\frac{t_{n,j}
+i\operatorname{sgn}(\iota_j(c))
\Delta_{n,j}}{2}
\right)^{-2}.
\end{equation}
Using $u=2\cos(\pi/5)$ and
$u'=2\cos(3\pi/5)$, we obtain
Table~\ref{tab:a5-local}.

\begin{table}[htbp]
\centering
\begin{tabular}{c|cc}
$n$&$d_{n,1}$&$d_{n,u}$\\\hline
$2$&$1$&$1$\\
$3$&$e^{2\pi i/3}$&$1$\\
$5$&$e^{2\pi i/5}$&$e^{4\pi i/5}$
\end{tabular}
\caption{Tangent determinants at the elliptic representatives.}
\label{tab:a5-local}
\end{table}

For $g_{n,c}$, the fixed cosets over
$\Gamma p_{n,c}$ are represented by
the identity and the nontrivial coset in
$N_Q(Q_{[p_{n,c}]})/Q_{[p_{n,c}]}$.
Conjugation by a representative of the
nontrivial coset fixes $g_{n,c}$ for
$n=2$ and sends it to $g_{n,c}^{-1}$
for $n=3,5$, since the centralizer of
an element of order $3$ or $5$ in $A_5$
is its cyclic subgroup.
Thus the tangent determinants of
$g_{n,c}$ at its two fixed points over
$\Gamma p_{n,c}$ are both $d_{n,c}$
for $n=2$, and are $d_{n,c}$ and
$d_{n,c}^{-1}$ for $n=3,5$.
These are also the determinants on $M$,
since $\HH^2\to M$ is locally
biholomorphic.

Within each $Q$-orbit, conjugation
identifies fixed points and conjugates
their tangent actions, preserving the
sum of local contributions.
Since $g_{n,1}$ and $g_{n,u}$ are
conjugate, we may compute the contributions
above $\Gamma p_{n,c}$ using $g_{n,c}$,
and add them to obtain $\chi_m(g_{n,1})$.

\begin{proposition}\label{prop:a5-nonidentity}
For $m\geq2$, the nonidentity character
values are
\begin{equation}\label{eq:a5-nonidentity}
\begin{aligned}
\chi_m(2A)&=4,\\
\chi_m(3A)&=
\begin{cases}
4,&3\mid m,\\
1,&3\nmid m,
\end{cases}\\
\chi_m(5A)=\chi_m(5B)&=
\begin{cases}
4,&5\mid m,\\
-1,&5\nmid m.
\end{cases}
\end{aligned}
\end{equation}
Here $5A$ and $5B$ are represented
by $g_{5,1}$ and $g_{5,1}^2$,
respectively.
\end{proposition}

\begin{proof}
The fixed-point formula and the
preceding calculation give
\[
\begin{aligned}
\chi_m(2A)&=4,\\
\chi_m(3A)&=2+e^{2\pi im/3}
                  +e^{-2\pi im/3},\\
\chi_m(5A)&=\sum_{r=1}^4 e^{2\pi irm/5}.
\end{aligned}
\]
These sums of roots of unity yield
the asserted values for $2A$, $3A$,
and $5A$.

The elements $g_{5,1}$ and $g_{5,1}^2$
belong to distinct conjugacy classes,
since their traces in
$\mathrm{SL}_2(\mathbb F_4)$
are $\bar u$ and $\bar u^2$.
They have the same fixed points,
and squaring the tangent actions
permutes the four nontrivial fifth
roots of unity. Hence
\[
\chi_m(5B)
=\sum_{r=1}^4e^{2\pi irm/5}
=\chi_m(5A).
\]
\end{proof}
\subsection{The equivariant index and its consequences}
\label{subsec:a5-results}

The character at the identity follows from the ordinary index
pairing. Together with the nonidentity values computed above,
it determines the equivariant index and its cuspidal interpretation.

\begin{lemma}\label{lem:a5-identity}
For every $m\geq 2$, the character of the equivariant
index at the identity satisfies
\begin{equation}\label{eq:a5-identity}
  \chi_m(e)=\chi(M)=4.
\end{equation}
\end{lemma}

\begin{proof}
We forget the $Q$-action and use the index pairing
\eqref{eq:a5-pairing}. Choose a compact cusp
truncation $W\subset M$ whose boundary consists
of sufficiently deep cusp level sets $q=T$, so
that $s=\partial\beta$ is nonzero outside
$\operatorname{int}W$ and
$\beta=\varepsilon q^2$ near $\partial W$.
By transversality for smooth bundle sections
\cite[Theorem~3.10]{CohenIgaNorbury},
the vector field $v=\nabla\beta$ admits an
arbitrarily small perturbation $v'$, supported
in $\operatorname{int}W$, with finitely many
nondegenerate zeros $p_1,\ldots,p_N$.
Set $s'=(v'^\flat)^{1,0}$.
The straight-line homotopy from $s$ to $s'$
is fixed outside a compact set, and hence
\[
  \Theta_{\beta,c}=\Theta_{s',c}
  \quad\text{in }K_c^0(M).
\]

Choose pairwise disjoint coordinate balls
$U_j\Subset\operatorname{int}W$ around the
zeros $p_j$, and put $U=\coprod_jU_j$.
The Koszul operator
\[
  c_{s'}=(s'\wedge+(s'\wedge)^*)^+
\]
satisfies
$c_{s'}^*c_{s'}=|s'|^2I$ and is therefore
invertible on $W\setminus U$.
Excision and additivity give
\[
  K^0(W,W\setminus U)
  \cong\bigoplus_jK_c^0(U_j).
\]
Writing $\Theta_j$ for the resulting local
classes, naturality of the Kasparov product
yields
\begin{equation}\label{eq:identity-local-sum}
  \chi_m(e)=
  \sum_j\Theta_j\widehat\otimes_{C_0(U_j)}
  [D_m^{0,*}]_{U_j}.
\end{equation}

Fix $p_j$ and choose a local unitary frame
$e_1,e_2$ of $T^{*1,0}M$. Write
$s'=a_1e_1+a_2e_2$ and set $F=(a_1,a_2)$.
In the induced exterior-algebra bases,
the even-to-odd Koszul map is
\[
  c_{s'}^+=
  \begin{pmatrix}
    a_1&-\overline{a_2}\\
    a_2&\overline{a_1}
  \end{pmatrix}.
\]
This is the complex two-dimensional Bott
matrix evaluated at $F$.
Since $p_j$ is nondegenerate, $F$ is
homotopic near $p_j$, through maps nonzero
on the boundary, to its invertible
real-linear differential $A_j=dF_{p_j}$.
Thus
\[
  \Theta_j=A_j^*\mathfrak b,
\]
where $\mathfrak b\in K_c^0(\mathbb C^2)$
is the class represented by the standard
Koszul matrix.

The twisting line is trivial over $U_j$,
and the restricted Dolbeault class has the
Euclidean principal-symbol model.
By the Bott--Dirac $KK$-equivalence
\cite[Theorems~2.7 and~2.10]{KasparovIndex},
the local pairing is
\[
  \Theta_j\widehat\otimes[D_m^{0,*}]_{U_j}
  =\deg(A_j)
   \bigl(\mathfrak b\widehat\otimes
   [D_{\mathbb C^2}^{0,*}]\bigr).
\]
With the total-degree grading used here,
the pairing in parentheses is $+1$.
Indeed, for $F(z)=(z_1,z_2)$, the
Dolbeault--Koszul complex has local
cohomology $\mathcal O/(z_1,z_2)
\cong\mathbb C$ concentrated in
holomorphic degree two.
Consequently,
\[
  \Theta_j\widehat\otimes[D_m^{0,*}]_{U_j}
  =\operatorname{sgn}\det_{\mathbb R}A_j.
\]

To identify this sign, use orthonormal
holomorphic coordinates at $p_j$.
For a real tangent vector
$v=\sum_{k=1}^2
(u_k\partial_{x_k}+w_k\partial_{y_k})$,
the identification $v\mapsto(v^\flat)^{1,0}$
is given at the point by
\[
  (v^\flat)^{1,0}
  =\frac12\sum_{k=1}^2
  (u_k-iw_k)\,dz_k.
\]
It reverses orientation in each complex
coordinate and hence preserves orientation
in complex dimension two. Therefore,
\[
  \operatorname{sgn}\det_{\mathbb R}A_j
  =\operatorname{ind}_{p_j}(v').
\]

Finally, $v'=v=2\varepsilon q\nabla q$
near $\partial W$. Since $W$ is truncated
at $q=T$, the outward unit normal is
$\nu=\nabla q/|\nabla q|$, and
\[
  \langle v',\nu\rangle
  =2\varepsilon q|\nabla q|>0.
\]
The Poincar\'e--Hopf theorem for compact
manifolds with boundary
\cite[Theorem~2.3.1]{RobbinSalamon}
now gives
\[
  \chi_m(e)
  =\sum_j\operatorname{ind}_{p_j}(v')
  =\chi(W)=\chi(M).
\]
Hirzebruch's calculation for the level-two
congruence quotient
\cite[\S5, p.~158]{Hirzebruch76} yields
\[
  \chi(M)
  =[\Gamma:\Gamma_0]\,2\zeta_F(-1)
  =60\cdot2\cdot\frac1{30}=4.
\]
\end{proof}

For the irreducible representations
$\mathbf1,\mathbf3,\mathbf3',\mathbf4,\mathbf5$ of $A_5$,
we use the character table
\begin{equation}\label{eq:a5-character-table}
\begin{array}{c|rrrrr}
&1A&2A&3A&5A&5B\\\hline
\text{class size}&1&15&20&12&12\\\hline
\mathbf1&1&1&1&1&1\\
\mathbf3&3&-1&0&u&u'\\
\mathbf3'&3&-1&0&u'&u\\
\mathbf4&4&0&1&-1&-1\\
\mathbf5&5&1&-1&0&0
\end{array}
\end{equation}
where $\mathbf3,\mathbf3'$ are distinguished by their
values on the class of $g_{5,1}$.

\begin{theorem}\label{thm:a5-index}
For $m\geq2$, the equivariant index $I_m\in R(A_5)$ is
\begin{equation}\label{eq:a5-decomposition}
I_m=
\begin{cases}
4\mathbf1,&15\mid m,\\
2\mathbf1-\mathbf3-\mathbf3'+2\mathbf4,
    &3\mid m,\ 5\nmid m,\\
3\mathbf1-\mathbf4+\mathbf5,
    &5\mid m,\ 3\nmid m,\\
\mathbf1-\mathbf3-\mathbf3'+\mathbf4+\mathbf5,
    &3\nmid m,\ 5\nmid m.
\end{cases}
\end{equation}
In particular, $I_{m+15}=I_m$ and $\dim I_m=4$.
\end{theorem}

\begin{proof}
Lemma~\ref{lem:a5-identity} and
Proposition~\ref{prop:a5-nonidentity} determine the character
on all five conjugacy classes. For an irreducible character
$\chi_\rho$, the coefficient of the corresponding representation
is
\[
\langle\chi_m,\chi_\rho\rangle
=\frac1{60}\sum_{[g]}|[g]|\chi_m(g)
                         \overline{\chi_\rho(g)}.
\]
Substitution from \eqref{eq:a5-character-table} gives
\eqref{eq:a5-decomposition}. Negative coefficients are allowed
because $I_m$ is a virtual representation.
\end{proof}

\begin{proposition}\label{cor:a5-invariant}
The index on $Q$-invariant sections is
\begin{equation}\label{eq:a5-average}
\begin{aligned}
\operatorname{ind}(B_m^+|_{L^2(M,E_m)^Q})
&=\frac1{60}\sum_{g\in Q}\chi_m(g)\\
&=\frac{4+15\cdot4+20\chi_m(3A)
       +12\chi_m(5A)+12\chi_m(5B)}{60}\\
&=1+\mathbf1_{3\mid m}+2\mathbf1_{5\mid m}.
\end{aligned}
\end{equation}
Here $\mathbf1_{d\mid m}$ is $1$ if $d\mid m$ and $0$
otherwise. This is also the index on $\Gamma$-invariant
sections on $\HH^2$.
\end{proposition}

\begin{proof}
The averaging projection $60^{-1}\sum_{g\in Q}g$ commutes
with $B_m$ and projects its graded kernels onto their
$Q$-invariant parts. Its supertrace gives
\eqref{eq:a5-average}. Moreover, $Q$-invariant
$\Gamma_0$-automorphic sections are precisely the
$\Gamma$-automorphic sections.
\end{proof}

To identify the cuspidal index with automorphic forms,
let $S_{k_1,k_2}(\Lambda)$ be the Hilbert cusp-form space
of classical weight $(k_1,k_2)$ for
$\Lambda=\Gamma_0$ or $\Gamma$.
The spaces for $\Gamma_0$ carry the natural $Q$-action.

\begin{proposition}\label{cor:a5-cusp-forms}
For $m\geq2$, the following identity holds in $R(A_5)$:
\begin{equation}\label{eq:a5-four-cusp-forms}
\begin{aligned}
I_m={}&[S_{2m,2m}(\Gamma_0)]
      -[S_{2m+2,2m}(\Gamma_0)]\\
     &-[S_{2m,2m+2}(\Gamma_0)]
      +[S_{2m+2,2m+2}(\Gamma_0)].
\end{aligned}
\end{equation}
Consequently,
\begin{equation}\label{eq:a5-full-group}
\begin{aligned}
&\dim S_{2m,2m}(\Gamma)
-\dim S_{2m+2,2m}(\Gamma)\\
&\quad-\dim S_{2m,2m+2}(\Gamma)
+\dim S_{2m+2,2m+2}(\Gamma)\\
&=1+\mathbf1_{3\mid m}+2\mathbf1_{5\mid m}.
\end{aligned}
\end{equation}
\end{proposition}

\begin{proof}
By \eqref{eq:a5-pairing},
$I_m=\Ind_Q(D_{m,\cusp}^+)$.

In holomorphic multidegree
$\boldsymbol p=(p_1,p_2)\in\{0,1\}^2$, the coefficient line is
\[
L_{\boldsymbol p}
=K_1^{m+p_1}\otimes K_2^{m+p_2}.
\]
Put $a_j=m+p_j$. For $q_j=1$,
a direct calculation from \eqref{eq:dirac-laplacian} gives
\[
D_j^2=\frac12T_j^*T_j+(a_j-1),
\qquad
T_j=2y_j\partial_{\bar z_j}-i(a_j-1),
\]
where the adjoint is taken with respect to the
hyperbolic measure in the scalar realization of
Section~\ref{subsec:invertibility}.
Since $D^2=D_1^2+D_2^2$, every sector with
$q_1+q_2>0$ satisfies
\[
\|Dv\|^2\geq(m-1)\|v\|^2
\]
for compactly supported smooth sections $v$.
Cutoff approximation on the complete manifold $M$
extends this estimate to the $L^2$ domain.
Hence all cuspidal harmonic forms of positive
antiholomorphic degree vanish.

In antiholomorphic degree zero, harmonic sections are
holomorphic sections of $L_{\boldsymbol p}$.
Their $L^2$ norms agree, up to a constant, with the
Petersson norms, and the cuspidal condition requires
vanishing constant terms at every cusp.
Thus these sections are precisely the Hilbert cusp
forms of classical weight
$(2m+2p_1,2m+2p_2)$.
This identification is $Q$-equivariant. Summing over the four multidegrees
with signs $(-1)^{p_1+p_2}$ gives
\eqref{eq:a5-four-cusp-forms}.
Finally, taking $Q$-invariants and applying
Proposition~\ref{cor:a5-invariant} proves
\eqref{eq:a5-full-group}.
\end{proof}

The dimensions of Hilbert cusp-form spaces are
classical objects of study in number theory,
notably through the dimension formulas of Shimizu
and subsequent developments.
The identity \eqref{eq:a5-full-group} provides
an index-theoretic relation among these dimensions,
with the periodic terms arising from the elliptic
fixed-point contributions.

More generally, \eqref{eq:a5-four-cusp-forms}
gives an identity in $R(A_5)$, retaining the
representation-theoretic information before taking
invariants. Thus our construction connects
the higher index in operator $K$-theory with
the arithmetic theory of Hilbert cusp forms.

\end{document}